\documentclass[10pt]{amsart}

\usepackage{amsmath,amsfonts, latexsym,graphicx, amssymb, amscd}

\newtheorem{theorem}{Theorem}[section]
\newtheorem{proposition}[theorem]{Proposition}
\newtheorem{lemma}[theorem]{Lemma}
\newtheorem{lemma-definition}[theorem]{Lemma and Definition}
\theoremstyle{definition}
\newtheorem{definition}[theorem]{Definition}

\newtheorem{corollary}[theorem]{Corollary}

\theoremstyle{remark}

\numberwithin{equation}{section}

\def\H{\mathbb H}

\title{Picard groups for line bundles with connection}
\date{}
\author{Helmut A. Hamm}
\address{Mathematisches Institut
Fachbereich Mathematik und Informatik,
Einsteinstrasse 62,
D-48149 M\"unster,
Germany }
 \email{hamm@uni-muenster.de}
\author{ L\^e D\~ung Tr\'ang }
\address{Universit\'e d'Aix-Marseille, 
LATP, UMR-CNRS 7353,
39 rue Joliot-Curie,
F-13453 Marseille Cedex 13, France}
 \email{ledt@ictp.it}

\begin{document}
\maketitle

\noindent{\bf Introduction}

\vskip.1in
Let $X$ be a complex analytic space. It is known that the isomorphism classes of
line bundles on $X$ define a group, called  the analytic Picard group $Pic^{an}(X)$ of $X$.

If  $X$ is a complex manifold, it is natural to consider 
line bundles on the space $X$ with connection or with integrable connection. The isomorphism classes
of these line bundles define 
groups that we shall denote by $Pic^{an}_{c}(X)$ for line bundles with connection
and $Pic^{an}_{ci}(X)$ for line bundles with integrable connection.

In this paper we are going to compare these groups with the original Picard group $Pic^{an}(X)$.

We first show that the groups $Pic^{an}_c(X)$ and $Pic^{an}_{ci}(X)$ are in general distinct from 
the Picard group.

We shall give a particular interest to the case of a non-singular complex algebraic variety $X$ for which we
can define the algebraic Picard group $Pic(X)$ of algebraic line bundles or the Picard group $Pic_{c}(X)$
(resp. $Pic_{ci}(X)$ or $Pic_{cir}(X)$) of algebraic line bundles
with connection (resp. integrable or regular integrable connection).

Note that $C^\infty$ connections are a standard tool in differential geometry. They are related to the differential-geometric description of characteristic classes - here it is important not to restrict to integrable connections. Complex analytic connections on principal bundles have been already studied by Atiyah \cite{A}, on vector bundles by Deligne \cite{[D]}. On Stein manifolds there are always complex analytic connections and these connections can also be used to define the complex 
first Chern class of a line bundle. 
The theory of $D$-modules which is only related to the integrable case will not 
be considered here. Finally there is also a relation to Deligne cohomology \cite{EV} which can be used to obtain results concerning the analytic case, but we have chosen a procedure which also carries directly over to the algebraic case.

The isomorphism classes of line bundle are in one-to-one correspondence with isomorphism classes of
invertible ${\mathcal O}_X$-modules (in the algebraic and analytic cases), this is why we shall use the 
isomorphism $Pic(X)\simeq H^1(X,{\mathcal O}^*_X)$ (resp. $Pic^{an}(X)\simeq H^1(X,{\mathcal O}^*_X)$).
In what follows we prefer mainly to use invertible sheaves instead of line bundles.

\section{Differentiable complex line bundles}

First, consider the differentiable case. Let $X$ be a paracompact differentiable (i.e. $C^\infty$) mani\-fold and 
let ${\mathcal E}^p_X$ be the sheaf of complex-valued differentiable $p$-forms on 
$X$, ${\mathcal E}_X:={\mathcal E}^0_X$. The isomorphism classes of differentiable complex 
line bundles form a group $Pic^\infty(X)\simeq H^1(X,{\mathcal E}^*_X)$, where ${\mathcal E}^*_X$ is 
the sheaf of non vanishing differentiable functions. 

Because of the exact exponential sequence $0\to \mathbb{Z}_X\to {\mathcal E}_X\stackrel{f\mapsto exp(2\pi if)}{\to} {\mathcal E}^*_X\to 0$ we obtain that the first Chern class gives an isomorphism $Pic^\infty(X)\simeq H^2(X;\mathbb{Z})$: notice that the
sheaf ${\mathcal E}_X$ is fine (see Exemple II 3.7.1 of \cite{Go}, p. 157f.), so $H^k(X,{\mathcal E}_X)=0$ for $k>0$, and the connecting homomorphism $H^1(X,\mathcal{E}_X^*)\to H^2(X;\mathbb{Z})$ may be used to define the first Chern class, see \cite{Hi} Theorem 4.3.1, p. 62.

Let $\mathcal{L}$ be an invertible $\mathcal{E}_X$-module which is the sheaf of sections of a line bundle $L$, $\mathcal{U}=(U_i)$ an open covering of $X$ such that $L|U_i$ is trivial. 

Let $g_{ij}$ be the transition function between $U_i$ and $U_j$ with respect to given trivializations $\phi_i:L|U_i\to U_i\times\mathbb{C}$, see \cite{Hi} p. 40 \S 3.2: $\phi_i\circ\phi_j^{-1}(x,t)=(x,g_{ij}(x)t)$, with $g_{ij}:U_i\cap U_j\to\mathbb{C}^*$. Let $s_i$: $s_i(x):=\phi_i^{-1}(x,1)$ be the corresponding section of $L|U_i$. Then we have 
$$s_j(x)=\phi_j^{-1}(x,1)=\phi_i^{-1}(x,g_{ij}(x))=g_{ij}(x)s_i(x)$$ 
i.e. $s_j=g_{ij}s_i$ (and not $s_i=g_{ij}s_j$, as suggested by \cite{GH} p. 70 !).\\

A connection on $\mathcal{L}$ is a $\mathbb{C}$-linear morphism $\nabla:{\mathcal L}\to {\mathcal E}^1_X\otimes_{{\mathcal E}_X} {\mathcal L}$ such that $\nabla(fs)=f\nabla(s)+df\otimes s$ (see \cite{GH} 0.5, p. 72). 
On $U_i$ it is given by $\alpha_i\in H^0(U_i,{\mathcal E}^1_X)$, where $\alpha_i$ is defined by 
$\nabla(s_i)=\alpha_i\otimes s_i$. 
We say then that $\nabla$ is represented by $(\alpha_i)$ with respect to $(s_i)$. 

We observe:

\begin{lemma}\label{crit}
There is a connection on the invertible $\mathcal{E}_X$-module $\mathcal{L}$ if and only if 
there are differential forms $\alpha_i$ defined on the open set $U_i$ such that on $U_i\cap U_j$ we have:
$$\frac{dg_{ij}}{g_{ij}}=\alpha_j-\alpha_i$$
\end{lemma}

\noindent {\bf Proof:} Let ($\alpha_i$) be a family of differential forms each one defined on the open space $U_i$.
They define on each ${\mathcal L}|U_i$ a connection $\nabla_i$:
$$\nabla_i(fs_i)=df\otimes s_i+f\alpha_i\otimes s_i.$$
Since $s_j=g_{ij}s_i=g_{ij}g_{ki}s_k=g_{kj}s_k$ on $U_i\cap U_j\cap U_k$ where $g_{ij}$ is
the transition function of ${\mathcal L}$ from $U_i$ to $U_j$ , we have 
$g_{ik}=g_{ij}g_{jk}$ on $U_i\cap U_j\cap U_k$. So the family
$(g_{ij})$ defines a 2-cocycle in $C^1(\mathcal{U},{\mathcal E}^*_X)$. Now on $U_i\cap U_j$,
we have $\nabla_i=\nabla_j$, for any $i,j$, if and only if, for any $i,j$: 
$$\nabla_i(s_j)=\nabla_i(g_{ij}s_i)=dg_{ij}\otimes s_i+g_{ij}\alpha_i\otimes s_i=\nabla_j(s_j)=\alpha_j\otimes s_j=g_{ij}\alpha_j\otimes s_i$$
or equivalently ($\alpha_i$) defines
a connection on ${\mathcal L}$ if and only if, for any $i,j$, on $U_i\cap U_i$ 
(compare with the proof of Theorem \ref{SE} below):
$$\frac{dg_{ij}}{g_{ij}}=\alpha_j-\alpha_i$$

The connection $\nabla=:\nabla^1$ induces a $\mathbb{C}$-linear mapping $\nabla^2:{\mathcal E}^1_X\otimes_{{\mathcal E}_X} {\mathcal L}\to {\mathcal E}^2_X\otimes_{{\mathcal E}_X} {\mathcal L}$, similarly as in the analytic case discussed in section 2.3. Then $\nabla^2\circ\nabla^1$ is ${\mathcal E}_X$-linear, it is the multiplication by a $2$-form in $H^0(X,{\mathcal E}^2_X)$, the curvature of 
$\nabla$. Its restriction to $U_i$ is $d\alpha_i$ (see \cite{GH} 0.5, p. 74/75). Then $\nabla$ is called integrable if and only if the curvature vanishes, i.e. $d\alpha_i=0$ for all $i$. 

The group of isomorphism classes of differentiable line bundles with an integrable connection is $Pic^\infty_{ci}(X)$. Let $c_1(\mathcal{L})_\mathbb{C}$ be the complex first Chern class of the line bundle $\mathcal{L}$.\\

Let $X$ be a paracompact differentiable manifold, we have:

\begin{lemma}\label{fund} a) There are connections on the invertible $\mathcal{E}_X$-module
$\mathcal{L}$  (see \cite{GH} 0.5, p. 73).\\
b) If $\nabla$ is a connection on $\mathcal{L}$, the curvature of $\nabla$ represents the element of: 
$$H^2_{DR}(X):=H^2(H^0(X,\mathcal{E}^\cdot_X))$$ 
which corresponds to $-2\pi ic_1(\mathcal{L})_\mathbb{C}\in H^2(X;\mathbb{C})$ under the De Rham isomorphism.\\
(Cf. \cite{GH} Proposition p. 141.)\\
c) The invertible $\mathcal{E}_X$-module $\mathcal{L}$ admits an integrable connection if and only if the complex first Chern class of $\mathcal{L}$ vanishes.\\
\end{lemma}

\noindent {\bf Proof:} a) The invertible sheaf $\mathcal{L}$ is given by a cocycle ($g_{ij})$, i.e. $s_j=g_{ij}s_i$, where $s_i\in H^0(U_i,\mathcal{L})$ is chosen as above. Then ($\frac{dg_{ij}}{g_{ij}}$) represents an element of $H^1({\mathcal U},{\mathcal E}^1_X)$. Now the \v{C}ech cohomology group $\check{H}^1(X,{\mathcal E}^1_X)\simeq H^1(X,{\mathcal E}^1_X)=0$, and $H^1({\mathcal U},{\mathcal E}^1_X)\to \check{H}^1(X,{\mathcal E}^1_X)$ is injective, see \cite{Fo} p. 91. So $H^1({\mathcal U},{\mathcal E}^1_X)=0$, and we can find $\alpha_i\in H^0(U_i,{\mathcal E}^1_X)$ such that $\alpha_j-\alpha_i=\frac{dg_{ij}}{g_{ij}}$, hence we get a connection on $\mathcal{L}$.\\
b) For the proof we hesitate to refer to \cite{GH} because of the difficulty in the relation between connections and transition functions mentioned above.\\ 
Anyhow the correctness of signs, however, is confirmed by \cite{BH} Theorem 29.4, p. 365, in connection with \cite{Hi} Th. 4.3.1, p. 62.\\
Let ${\mathcal U}=(U_i)_i$ be a simple covering (``recouvrement simple"), i.e. an open covering of $X$ such that the finite intersections of the $U_i$ are empty or contractible, see \cite{W} p. 120 (for the notion of contractibility cf. footnote on \cite{BT} p. 36). Then $\mathcal{L}|U_i$ is trivial for all $i$, cf, \cite{BT} Cor. 6.9, p. 59. 
From Lemma \ref{crit} above, we have  
$\frac{dg_{ij}}{g_{ij}}=\alpha_j-\alpha_i$, where ($\alpha_i$) represents the connection. 

Following 
\cite{GH} 1.1, p. 141 adapted to complex valued forms, ($d\alpha_i$) defines the curvature in $H^0(X,d{\mathcal E}^1_X)$. It defines an element in: 
$$H^2_{DR}(X)=H^0(X,d{\mathcal E}^1_X)/dH^0(X,{\mathcal E}^1_X)\simeq H^0({\mathcal U},d{\mathcal E}^1_X)/dH^0({\mathcal U},{\mathcal E}^1_X).$$
We have another way to consider the De Rham cohomology.\\
Look at the double complex $C^\cdot(\mathcal{U},\mathcal{E}_X^\cdot)$.\\
Let $(C^\cdot(\mathcal{U},\mathcal{E}_X^\cdot))_{tot}$ be the corresponding total (or: simple) complex: 
$$(C^\cdot(\mathcal{U},\mathcal{E}_X^\cdot))^k_{tot}:=\oplus_{p+q=k} C^p(\mathcal{U},\mathcal{E}_X^q)$$
The coboundary operator $d^+$ is defined as follows, see \cite{Go} II 4.6, p. 176: 
$$d^+|C^p(\mathcal{U},\mathcal{E}_X^q)=\delta+(-1)^pd,$$ 
where
$\delta:C^p(\mathcal{U},\mathcal{E}_X^q)\to C^{p+1}(\mathcal{U},\mathcal{E}_X^q)$ is the \v{C}ech coboundary and 
$d:C^p(\mathcal{U},\mathcal{E}_X^q)\to C^p(\mathcal{U},\mathcal{E}_X^{q+1})$ the coboundary operator of the complex of differential forms.\\
We have:
\begin{lemma}
$H^k_{DR}(X)\simeq H^k(H^0(\mathcal{U},\mathcal{E}_X^\cdot))
\simeq H^k(C^\cdot(\mathcal{U},\mathcal{E}_X^\cdot)_{tot})$.
\end{lemma}
\noindent {\bf Proof:} We notice that we have a spectral sequence:
$$H^k(H^l(\mathcal{U},\mathcal{E}_X^\cdot))
\Rightarrow H^{k+l}(C^\cdot(\mathcal{U},\mathcal{E}_X^\cdot)_{tot}).$$
Th\'eor\`eme 5.2.3 b) of \cite{Go} implies that $H^l(\mathcal{U},\mathcal{E}_X^\cdot)=0$, for $l>0$
because the sheaves ${\mathcal E}^k$ are fine.

Therefore $H^k(H^0(\mathcal{U},\mathcal{E}_X^\cdot))=H^k(C^\cdot(\mathcal{U},\mathcal{E}_X^\cdot)_{tot})$.

Using again \cite{Go} but Th\'eor\`eme 5.2.2, we have $H^0(\mathcal{U},\mathcal{E}_X^p)
=H^0(X,{\mathcal E}^p)$. So $H^0(\mathcal{U},\mathcal{E}_X^\cdot)$ is the De Rham complex and
$H^k_{DR}(X)\simeq H^k(H^0(\mathcal{U},\mathcal{E}_X^\cdot))$ 

\vskip.1in
So, the class of the curvature corresponds to the element of $H^2(C^\cdot(\mathcal{U},\mathcal{E}_X^\cdot)_{tot})$ represented by $(d\alpha_i)\in C^0(\mathcal{U},\mathcal{E}_X^2)$.\\
Looking at the coboundary of $(\alpha_i)\in C^0(\mathcal{U},\mathcal{E}_X^1)$, viewed as element of 
$(C^\cdot(\mathcal{U},\mathcal{E}_X^\cdot))_{tot}^1$, we see that $(d\alpha_i)$ is 
cohomologous to $-(\alpha_j-\alpha_i)=-(\frac{dg_{ij}}{g_{ij}})$.\\
Looking at the coboundary of $(log\,g_{ij})\in C^1(\mathcal{U},\mathcal{E}^0_X)$ in the total complex (here $log\,g_{ij}$ denotes some branch of the logarithm of $g_{ij}$, i. e. a suitable antiderivative of $\frac{dg_{ij}}{g_{ij}}$ ,
see Theorem 4.3.1 p. 62 of \cite{Hi}), we see that $-(\frac{dg_{ij}}{g_{ij}})$ is cohomologous to $-(log\,g_{jk}-log\,g_{ik}+log\,g_{ij})$.\\
Since this element is a cocycle in $(C^\cdot(\mathcal{U},\mathcal{E}_X^\cdot))_{tot}$, its differential must be $0$,
so we have that $c_{ij}=(log\,g_{jk}-log\,g_{ik}+log\,g_{ij})\in\mathbb{C}$.\\
The class of the curvature of ${\mathcal L}$ is therefore equal to the class given by the constants 
$-c_{ij}=-(log\,g_{jk}-log\,g_{ik}+log\,g_{ij}).$\\
On the other hand, $c_1(\mathcal{L})$ may be defined as the image of the class 
$[\mathcal{L}]\in H^1(X,\mathcal{O}^*_X)$ in $H^2(X,\mathbb{Z}_X)$ with respect to the exponential sequence.
Since, by \cite{W}, we can choose a simple covering ${\mathcal U}$ of $X$, we have an exact sequence:
$$0\to C^\cdot ({\mathcal U},{\mathbb Z})\to C^\cdot ({\mathcal U},{\mathcal O}_X)\to 
C^\cdot ({\mathcal U},{\mathcal O}^*)\to 0$$ 
Therefore we obtain a homomorphism $H^1(\mathcal{U},\mathcal{O}^*_X)\to H^2(\mathcal{U},\mathbb{Z}_X))$. The image of the class of ${\mathcal L}$, given by the transition functions $g_{ij}$ on the covering ${\mathcal U}$, in $H^1(\mathcal{U},\mathcal{O}^*_X)$ in 
$H^2(\mathcal{U},\mathbb{Z}_X))$ is $\frac{1}{2\pi i}(log\,g_{jk}-log\,g_{ik}+log\,g_{jk})$. \\
The class of the curvature of $\nabla$ in $H^2_{DR}(X)$ corresponds under De Rham isomorphism to
$-2\pi ic_1({\mathcal L})_{\mathbb C}$.\\
c) Note that $Pic^\infty_{ci}(X)\simeq H^1(X;\mathbb{C}^*)$ (compare with Proposition \ref{iso} below).
 Then we have a commutative diagram:
$$\begin{array}{ccccc}
H^1(X;\mathbb{C}^*)&\to&H^2(X;\mathbb{Z})&\to&H^2(X;\mathbb{C})\\
\downarrow&&\downarrow=&&\\
H^1(X,{\mathcal E}^*_X)&\stackrel{\simeq}{\to}&H^2(X;\mathbb{Z})&&
\end{array}$$
where the upper row is exact; for the isomorphism in the lower row see beginning of this section. This implies our statement.\\

\section{Analytic Comparisons}

\subsection{}{\bf $Pic^{an}(X)$ and $Pic^{an}_c(X)$}\label{GR} \\

In this section let $X$ be a complex manifold which is paracompact (e.g. Stein or compactifiable; the condition is not automatically fulfilled, see \cite {CR}). Similarly as before, a connection on an invertible $\mathcal{O}_X$-module $\mathcal{L}$ is a $\mathbb{C}$-linear morphism $\nabla:{\mathcal L}\to \Omega^1_X\otimes_{{\mathcal O}_X} {\mathcal L}$ such that $\nabla(fs)=f\nabla(s)+df\otimes s$, see \cite{[D]} I D\'ef. 2.4, p. 7. We obtain groups $Pic^{an}_c(X)$.\\

We have an exact sequence of sheaves:
$$0\to {\mathbb C}^*_X\to {\mathcal O}_X^*\to d{\mathcal O}_X\to 0$$
where ${\mathbb C}^*_X\to {\mathcal O}_X^*$ is given by the inclusion and ${\mathcal O}_X^*\to d{\mathcal O}_X$ is defined by $f\mapsto df/f$.

This latter morphism is surjective, because, if $\omega\in d{\mathcal O}_{X,x}$, there 
is $f\in{\mathcal O}_{X,x}$ such that $\omega=df$. Then $e^f\in{\mathcal O}_{X,x}^*$ has its image equal to 
$\omega$. The rest of the sequence is exact because of Poincar\'e Lemma.

This exact sequence of sheaves gives an exact sequence of cohomology:
$$\ldots \to H^0(X,d{\mathcal O}_X)\to H^1(X,{\mathbb C}^*_X)\to H^1(X,{\mathcal O}^*_X)\to
H^1(X,d{\mathcal O}_X)\to \ldots$$
Here we only use the mapping 
$H^1(X,{\mathcal O}^*_X)\to H^1(X,d{\mathcal O}_X)$ from
the exact sequence (see also the proof of Theorem 2.2.22 of \cite{Br}):

\begin{theorem} \label{SE} We have an exact sequence:
$$H^0(X,\mathcal{O}^*_X)\to H^0(X,\Omega^1_X)\to Pic^{an}_{c}(X)\to Pic^{an}(X)\to H^1(X,\Omega^1_X)$$
\end{theorem}

\noindent{\bf Proof}.
The map $H^0(X,\mathcal{O}^*_X)\to H^0(X,\Omega^1_X)$ is defined by $g\mapsto\frac{dg}{g}$, $H^0(X,\Omega^1_X)\to Pic^{an}_{c}(X)$ by 
$\omega\mapsto ({\mathcal O}_X, \nabla(f)=df +f\omega)$.

The map $Pic^{an}(X)\to H^1(X,\Omega^1_X)$ is the composition of 
$H^1(X,{\mathcal O}^*_X)\to
H^1(X,d{\mathcal O}_X)$ and the natural map from $H^1(X,d{\mathcal O}_X)$ to
$H^1(X,\Omega^1_X)$, since $Pic^{an}(X)\simeq H^1(X,{\mathcal O}^*_X)$.

Now, notice that we have a group structure on $Pic^{an}_{c}(X)$. According to Deligne
in  \cite{[D]} p. 8, 
consider the invertible sheaves (i.e. invertible $\mathcal{O}_X$-modules) ${\mathcal L}$
and ${\mathcal L}'$ defined by the ($s_i$) and ($s_i'$) on an open covering ${\mathcal U}$, 
with the connections $\nabla$ and $\nabla'$ defined by 
($\alpha_i$) and ($\alpha_i'$) on ${\mathcal U}$, then
${\mathcal L}\otimes{\mathcal L}'$ is invertible and defined by ($s_i\otimes s_i'$), and the connection $\nabla_0$ on this invertible sheaf is defined by ($\alpha_i+\alpha_i'$).

Now let us prove the exactness. First, $g\in H^0(X,\mathcal{O}^*_X)$ is mapped to $\frac{dg}{g}\in H^0(X,\Omega^1_X)$, and this in turn to the element of $Pic^{an}_c(X)$ represented by $(\mathcal{O}_X,\nabla)$, where $\nabla(f):=df+\frac{dg}{g}f$. This is the inverse image of $(\mathcal{O}_X,d)$ under the isomorphism $\cdot g:\mathcal{O}_X\to \mathcal{O}_X$, so its class in $Pic^{an}_c(X)$ is trivial:
$$\begin{array}{cccc}&{\mathcal O}_X&\stackrel{\nabla}\to&\Omega_X^1=\Omega_X^1\otimes {\mathcal O}_X\\
&{.g}\downarrow&&{.g}\downarrow\\
&{\mathcal O}_X&\stackrel{d}\to&\Omega_X^1=\Omega_X^1\otimes {\mathcal O}_X\\
\end{array}$$
Suppose now that $\omega\in H^0(X,\Omega^1_X)$ is mapped onto the trivial element of $Pic_c^{an}(X)$, which means that $(\mathcal{O}_X,d)$ is isomorphic to $(\mathcal{O}_X,f\mapsto df+\omega f)$. The isomorphism is of the form $\cdot g$ for some $g\in H^0(X,\mathcal{O}^*_X)$. Then, the image of $1$ is $\omega\in \Omega^1_X$ and by the multiplication by $g$, it is $dg$. Therefore $\omega=\frac{dg}{g}$.

The kernel of the map $Pic^{an}_{c}(X)\to Pic^{an}(X)$ defined by 
$({\mathcal L},\nabla)\mapsto {\mathcal L}$ is given by the pairs $({\mathcal O}_X,\nabla)$, so 
it coincides with the image of the morphism $H^0(X,\Omega^1_X)\to Pic^{an}_{c}(X)$ defined by 
$\omega\mapsto ({\mathcal O}_X, \nabla(f)=df +f\omega$). So, the middle part of the sequence is exact.

Now let $\mathcal L$ be an invertible sheaf which is in the kernel of $Pic^{an}(X)\to H^1(X,\Omega^1_X)$.
Let ${\mathcal U}=(U_i)_{i\in I}$ be a covering of $X$, such that ${\mathcal L}|U_i$ is isomorphic to
${\mathcal O}|U_i$ by a map  ${\mathcal O}_X|U_i\rightarrow {\mathcal L}|U_i$ which 
corresponds to $1\mapsto s_i$. Let $g_{ij}$ be the complex analytic transition map defined on $U_{ij}=U_i\cap U_j$ from ${\mathcal L}|U_i$ to ${\mathcal L}|U_j$. 
We have $s_j=g_{ij}s_i$ on $U_i\cap U_j$. 

Since $s_j=g_{ij}s_i=g_{ij}g_{ki}s_k=g_{kj}s_k$ on $U_i\cap U_j\cap U_k$, we have 
$g_{kj}=g_{ij}g_{ki}$ on $U_i\cap U_j\cap U_k$. The family
$(g_{ij})$ defines a 2-cocycle of $H^1(X,{\mathcal O}^*_X)$, a fact which is well-known. 
Since $H^1({\mathcal U},\Omega^1_X)\subset H^1(X,\Omega^1_X)$, the image of ${\mathcal L}$
in $H^1(X,\Omega^1_X)$ being trivial, the $2$-cocycle $(dg_{ij}/g_{ij})$ is trivial, i.e. a coboundary. 
Therefore 
there are differential forms $\omega_i$ and $\omega_j$ defined respectively 
on $U_i$ and $U_j$, such that:
$${dg_{ij}\over g_{ij}}=\omega_j-\omega_i$$
on $U_i\cap U_j$.

Consider for each $i$ the connection $\tilde{\nabla}_i$ on ${\mathcal O}_X|U_i$ defined by:
$$\tilde{\nabla}_i(f)=df+f\omega_i$$
This defines on ${\mathcal L}|U_i$ a connection:
$$\nabla_i(fs_i)=df\otimes s_i+f\omega_i\otimes s_i,$$
which gives for $f=1$:
$$\nabla_i(s_i)=\omega_i\otimes s_i.$$

On $U_i\cap U_j$, we have $ g_{ij} s_i= s_j$. Therefore, on $U_i\cap U_j$:
$$\nabla_i(fg_{ij} s_i)=d(fg_{ij})\otimes s_i+fg_{ij}\omega_i\otimes s_i=g_{ij}df\otimes s_i+fdg_{ij}\otimes s_i+fg_{ij}\omega_i\otimes s_i,$$
which implies, with $f=1$, on $U_i\cap U_j$:
$$\nabla_i(g_{ij} s_i)=dg_{ij}\otimes s_i+g_{ij}\omega_i\otimes s_i.$$
Therefore:
$$\nabla_i(s_j)=g_{ij}({dg_{ij}\over g_{ij}}+\omega_i)\otimes s_i=(\omega_j-\omega_i+\omega_i)\otimes s_j$$
which yields:
$$\nabla_i(s_j)=\nabla_j(s_j)$$
on $U_i\cap U_j$.

Therefore the $(\nabla_i)_{i\in I}$ define on ${\mathcal L}$ a connection $\nabla$ and the 
class of the element 
${\mathcal L}$
which lies in the kernel of the map $Pic^{an}(X)\to H^1(X,\Omega^1_X)$ 
is the image of the class of $({\mathcal L},\nabla)$.

\vskip.1in
It remains to prove that the image of $({\mathcal L},\nabla)$ in $H^1(X,\Omega^1_X)$ in the above sequence vanishes.

\vskip.1in
Let $(U_i)_{i\in I}$ be an open covering of $X$ such that ${\mathcal L}|U_i$ is isomorphic
to ${\mathcal O}_X|U_i$ by a map $s_i\mapsto 1$. We write $\nabla s_i=\omega_i\otimes s_i$. Let $(g_{ij})$ be the cocycle of transition functions
such that $s_j=g_{ij}s_i$. Then $(dg_{ij}/g_{ij})$ is a cocycle which represents an element of $H^1(X,\Omega^1_X)$. Since: 
$$\nabla(s_j)=\nabla(g_{ij}\otimes s_i)=dg_{ij}\otimes s_i+g_{ij}\omega_i\otimes s_i=\omega_j\otimes s_j
=g_{ij}\omega_j\otimes s_i,$$
we obtain:
$${dg_{ij}\over g_{ij}}=\omega_j-\omega_i.$$
Therefore the class of the element given by the elements $(dg_{ij}/ g_{ij})$ vanishes in $H^1(X,\Omega^1_X)$.

This shows that the above sequence is exact.

\vskip.1in 
We shall give an interpretation of this exact sequence below.

\vskip.1in
Implicitly we have used:

\begin{lemma} \label{Cech}
Let ${\mathcal L}$ be an invertible $\mathcal{O}_X$-module which is represented by a cocycle $(g_{ij})$
in $C^1({\mathcal U},{\mathcal O}^*_X)$. Then, the connection $\nabla$ on ${\mathcal L}$ is represented by an element $(\omega_i)$ in $C^0({\mathcal U}, \Omega^1_X)$ which is mapped by 
$\delta:C^0({\mathcal U}, \Omega^1_X)\to C^1({\mathcal U},\Omega^1_X)$ 
onto $(\frac{dg_{ij}}{g_{ij}})\in C^1({\mathcal U},\Omega^1_X)$.
\end{lemma}

Note that $(d\omega_i)\in {\mathcal C}^0({\mathcal U}, \Omega^2_X)$ is a cocycle, i.e. 
an element of $H^0(X,\Omega^2_X)$, which is the curvature of $\nabla$, see below.\\

Particularly easy is the case of Stein manifolds. Using the same methods as in the preceding section 
(Lemma \ref{fund}) we obtain:

\begin{lemma}\label{Stein}
Let $\mathcal{L}$ be an invertible $\mathcal{O}_X$-module on a Stein manifold $X$.\\
a) There is a complex analytic connection on $\mathcal{L}$.\\
b) There is an integrable complex analytic connection on $X$ if and only if the complex first Chern class of $\mathcal{L}$ vanishes.
\end{lemma}

In the following subsection we shall show how our reasoning above is related to the literature.

\subsection{}{\bf Differentiable connections and Atiyah obstruction.} \\

Atiyah (\cite{A} \S 2) has studied complex analytic connections on a holomorphic principal fibre bundle $P$. Whereas differentiable connections always exist there is an obstruction to the existence of a 
complex analytic one. In particular, there is an obstruction $b(E)$ to the existence of a complex analytic connection on the principal fibre bundle which corresponds to a holomorphic vector bundle $E$ (see \cite{A} p. 194). We call it the Atiyah obstruction. 
In the case of a line bundle $L$ we have that 
$b(L)\in H^1(X,\Omega^1_X)$. 

Here we use again invertible sheaves $\mathcal{L}$ instead of line bundles $L$. Then a complex analytic connection on $L$ corresponds to a connection on the sheaf $\mathcal{L}$ of holomorphic sections of $L$.

First we work in the analytic context. Let us recall the definition of $b(\mathcal{L})$, see \cite{A} p. 193. Let $D(\mathcal{L})$ be the locally free $\mathcal{O}_X$-module defined as follows:\\
as a $\mathbb{C}_X$-module, $D(\mathcal{L}):=\mathcal{L}\oplus(\Omega^1_X\otimes_{\mathcal{O}_X}\mathcal{L})$, and the ${\mathcal O}_X$-module structure is given by
$f\cdot(s,\beta):=(fs,f\beta+ df\otimes s)$, if $f$ is a section of ${\mathcal O}_X$, $s$ a section of ${\mathcal L}$
and $\beta$ is a section of $\Omega^1_X\otimes_{\mathcal{O}_X}\mathcal{L}$.\\
Then we get an exact sequence of $\mathcal{O}_X$-modules
$$0\to \Omega^1_X\otimes_{\mathcal{O}_X}\mathcal{L}\to D(\mathcal{L})\to\mathcal{L}\to 0$$
where the second arrow is given by $\beta\mapsto(0,\beta)$ and the third one by $(s,\beta)\mapsto s$.\\
Applying $Hom(\mathcal{L},\cdots)$ we obtain a long exact cohomology sequence
$$\ldots\to H^0(X,Hom(\mathcal{L},D(\mathcal{L}))\to H^0(X,Hom(\mathcal{L},\mathcal{L}))\to H^1(X,Hom(\mathcal{L},\Omega^1_X\otimes_{\mathcal{O}_X}\mathcal{L}))\to\ldots$$
Now $b(\mathcal{L})$ is defined as the image of $1\in H^0(X,\mathcal{O}_X)$ in $H^1(X,\Omega^1_X)$ with 
respect to the identifications
$H^0(X,\mathcal{O}_X)\simeq H^0(X,Hom(\mathcal{L},\mathcal{L}))$ and  $H^1(X,\Omega^1_X)\simeq H^1(X,Hom(\mathcal{L},\Omega^1_X\otimes_{\mathcal{O}_X}\mathcal{L}))$
(so the mapping depends on $\mathcal{L}$ !).

\begin{lemma}\label{b(L)} $b(\mathcal{L})=0$ if and only if $\mathcal{L}$ admits a connection.
\end{lemma}

\noindent{\bf Proof:} A splitting of the first exact sequence above is given by $\mathcal{L}\to D(\mathcal{L}): s\mapsto(s,\nabla(s))$ for some connection $\nabla$ on $\mathcal{L}$. 

This defines an element of $H^0(X,Hom(\mathcal{L},D(\mathcal{L}))$
whose image in $H^0(X,Hom(\mathcal{L},\mathcal{L}))$ is $1$ and by the exactness of the second sequence
its image in :
$$H^1(X,Hom(\mathcal{L},\Omega^1_X\otimes_{\mathcal{O}_X}\mathcal{L}))$$ 
is $0$.\\
On the other hand, a splitting corresponds to an inverse image of: 
$$1\in H^0(X,\mathcal{O}_X)\simeq H^0(X,Hom(\mathcal{L},\mathcal{L}))$$ 
in the second exact sequence. This proves the converse.

\begin{lemma}\label{atiyah} $b(\mathcal{L})$ is the image of $-[\mathcal{L}]\in H^1(X,\mathcal{O}^*_X)$ in $H^1(X,\Omega^1_X)$, i.e. $b(\mathcal{L})$ is represented by $-(\frac{dg_{ij}}{g_{ij}})$.
\end{lemma}

\noindent{\bf Proof:} Let ${\mathcal U}=(U_i)$ be an open Stein covering of $X$ such that ${\mathcal L}|U_i$ is trivial. Let $s_i$ be a nowhere vanishing section of ${\mathcal L}|U_i$. Then, $s_j=g_{ij}s_i$, 
where $g_{ij}$ are the corresponding transition functions. Let $\nabla_i$ be the connection on $\mathcal{L}|U_i$
such that $\nabla_i(s_i)=0$. Now, let us describe 
$H^0({\mathcal U},Hom({\mathcal L},{\mathcal L}))
\to H^1({\mathcal U},Hom({\mathcal L},\Omega^1_X\otimes_{\mathcal{O}_X}\mathcal{L}))$ using the exact sequence of complexes:
$$0 \to C^\cdot ({\mathcal U},Hom({\mathcal L},\Omega^1_X\otimes_{\mathcal{O}_X}\mathcal{L}))
\to C^\cdot({\mathcal U},Hom({\mathcal L},D({\mathcal L}))\to 
C^\cdot({\mathcal U},Hom({\mathcal L},{\mathcal L}))\to  0.$$
Consider $(\sigma_i)\in C^0({\mathcal U},Hom({\mathcal L},D({\mathcal L}))$, where $\sigma_i$ is the homomorphism ${\mathcal L}|U_i\to D({\mathcal L})|U_i$ which maps $s_i$ to $(s_i,0)$ 
(note that $\nabla_i(s_i)=0$), i.e. $s_j=g_{ij}s_i$ to $(s_j,\frac{dg_{ij}}{g_{ij}}\otimes s_j)$.
Then $(\sigma_i)$ is mapped to $(\tau_i)\in C^0({\mathcal U},Hom({\mathcal L},{\mathcal L}))$ 
with $\tau_i=id:{\mathcal L}|U_i\to{\mathcal L|U_i}$.

The coboundary of $(\sigma_i)$ is given by 
$\sigma_j-\sigma_i:{\mathcal L}|U_i\cap U_j\to D({\mathcal L})|U_i\cap U_j$:
$$(\sigma_j-\sigma_i)(s_j)=(0,-\frac{dg_{ij}}{g_{ij}}\otimes s_j),$$
so $\sigma_j-\sigma_i$ can be identified with $-\frac{dg_{ij}}{g_{ij}}\in H^0(U_i\cap U_j,\Omega^1_X)$.

\begin{corollary}\label{existence}
An invertible sheaf ${\mathcal L}$ admits a connection if and only if its image in $H^1(X,\Omega^1_X)$
is $0$.
\end{corollary}

This corollary is consequence of Lemmas \ref{b(L)} and \ref{atiyah}.
This coincides with our result from Theorem \ref{SE}.\\

Now let us turn to the differentiable context. Note that a differentiable connection on the holomorphic line bundle 
$L$ corresponds to a connection on the sheaf 
$\mathcal{E}_X\otimes_{\mathcal{O}_X}\mathcal{L}$ of differentiable sections of $L$.

We shall give another proof of Lemma \ref{atiyah}.

Precisely, in the case of invertible $\mathcal{O}_X$-modules it is natural to look at differentiable connections, i.e. connections on $\mathcal{E}_X\otimes_{\mathcal{O}_X}\mathcal{L}$, which are compatible with the complex structure (see \cite{GH} \S  0.5, p. 73) in the sense that 
$\nabla\,s-\bar{\partial}s\in {\mathcal E}^{1,0}_X\otimes_{{\mathcal O}_X} {\mathcal L}$. 
Locally a compatible connection is given by $\alpha_i\in H^0(U_i,{\mathcal E}^{1,0}_X)$, with respect to holomorphic sections $s_i\in H^0(U_i,\mathcal{L})$. Note that a complex connection (i.e. on $\mathcal{L}$) 
extends to such a differentiable one (i.e. on $\mathcal{E}_X\otimes_{\mathcal{O}_X}\mathcal{L}$) which is compatible with the complex structure:\\
For $f\in H^0(U_i,\mathcal{E}_X)$, put $\nabla(fs_i):=df\otimes s_i+f\nabla s_i$. In fact, $\nabla(fs_i)-\bar{\partial}(fs_i)=\partial f\otimes s_i+f\nabla s_i$ is an element of 
$H^0(U_i,\mathcal{E}^{1,0}_X\otimes_{{\mathcal O}_X}\mathcal{L})$.\\

\begin{lemma} [cf. \cite{GH} \S  0.5, p. 73] There are differentiable connections on $\mathcal{L}$ which 
are compatible with the complex structure.
\end{lemma}

\noindent {\bf Proof:} Let $g_{ij}$ be the transition functions which define $\mathcal{L}$. 
Since the sheaf ${\mathcal E}^{1,0}_X$
is an ${\mathcal E}_X$-module, it is a fine sheaf (see Exemple 3.7.1 of \cite{Go}),
so $H^1(X,{\mathcal E}^{1,0}_X)=0$. 
So we have that 
($\frac{dg_{ij}}{g_{ij}}$)$\in H^1(X,{\mathcal E}^{1,0}_X)$ is $0$, hence ($\frac{dg_{ij}}{g_{ij}}$)$\in H^1(\mathcal{U},{\mathcal E}^{1,0}_X)$, too. Therefore there is 
($\alpha_i$)$\in H^0(X,{\mathcal E}^{1,0}_X)$ such that the ($\frac{dg_{ij}}{g_{ij}}$) 
is a coboundary of ($\alpha_i$).
As we saw above, the connection given by ($\alpha_i$) must be compatible with the complex structure.

Alternative: Choose a Hermitian metric on the associated line bundle $L$ and take the corresponding metric connection, see Lemma \ref{curvmet} below.\\

Let ${\mathcal U}$ be an open Stein covering ($U_i$) of the complex manifold $X$. Let $(\alpha_i)\in C^0({\mathcal U},{\mathcal E}^{1,0}_X)$ represent a (differentiable) connection 
$\nabla$ relatively to the sections $s_i$ of
${\mathcal L}\subset \mathcal{E}_X\otimes_{\mathcal{O}_X}\mathcal{L}$ on $U_i$ which is compatible with the complex structure, where 
($g_{ij}$) is a cocycle defining $\mathcal{L}$ and: 
$$\alpha_j-\alpha_i=\frac{dg_{ij}}{g_{ij}}$$ 
on $U_i\cap U_j$. 
Let $\kappa$ be the $(1,1)$-part of the curvature of $\nabla$, so $\kappa$ coincides with $\bar{\partial}\alpha_i$ on $U_i$, the $\bar{\partial}$-curvature of $\nabla$. 

Then we have:

\begin{lemma}\label{curv} Let ${\mathcal L}$ be a complex analytic invertible sheaf and let 
$\nabla$ be a (differentiable) connection on $\mathcal{E}_X\otimes_{\mathcal{O}_X}\mathcal{L}$ 
which is compatible with the complex structure. Let $\kappa$ be the $\bar{\partial}$-curvature of $\nabla$. Then:
\begin{enumerate}
\item{} $\nabla$ ``is'' a complex analytic connection (i.e. $\nabla|\mathcal{L}$ is a connection) 
if and only if $\kappa=0$.
\item{} $b(\mathcal{L})$ is the Dolbeault cohomology class of the $\bar{\partial}$-curvature $\kappa$
of $\nabla$ in $H^1(X,\Omega^1_X)$.
\item{} $\mathcal{L}$ admits a complex analytic connection $\nabla'$ if and only if the cohomology class of $\kappa$ vanishes.
\end{enumerate}
\end{lemma}

\noindent {\bf Proof:} (1) $\kappa=0$ means that the $\alpha_i$ are holomorphic.\\
(2) Let $\nabla$ be given by $(\alpha_i)_i$. The Dolbeault cohomology class of $\kappa$ in 
$H^1(H^0(X,{\mathcal E}^{1,*}_X))\simeq H^0({\mathcal U},\bar{\partial}{\mathcal E}^{1,0}_X)/\bar{\partial}H^0({\mathcal U},{\mathcal E}^{1,0}_X)$ is represented by $(\bar{\partial}\alpha_i)$.\\
We may consider $C^0({\mathcal U},{\mathcal E}^{1,0})$ as the term of degree $0$ of the complex 
$(C^\cdot({\mathcal U},{\mathcal E}^{1,\cdot}))_{tot}$. Here, the coboundary operator $d^+$ is given as follows: $d^+|C^p(\mathcal{U},\mathcal{E}^{1,q})=\delta+(-1)^p\bar{\partial}$. The image of $(\alpha_i)$ in the term of degree $1$
of $C^\cdot({\mathcal U},{\mathcal E}^{1,\cdot})_{tot}$ is $((\alpha_j-\alpha_i)+\bar\partial(\alpha_i))$
in $C^1({\mathcal U},{\mathcal E}^{1,0})\oplus C^0({\mathcal U},{\mathcal E}^{1,1})$. This image is a coboundary
and defines the element $0$ in $H^1(C^\cdot({\mathcal U},{\mathcal E}^{1,\cdot})_{tot})$. The elements 
$(\alpha_j-\alpha_i)=(\frac{dg_{ij}}{g_{ij}})$ and $\bar\partial \alpha_i$ are cocycles in 
$C^1({\mathcal U},{\mathcal E}^{1,0})\oplus C^0({\mathcal U},{\mathcal E}^{1,1})$.
Therefore, the image of
$-(\frac{dg_{ij}}{g_{ij}})$ in $H^1((C^\cdot({\mathcal U},{\mathcal E}^{1,\cdot})_{tot})$ equals the
image of $\bar\partial \alpha_i$ in $H^1((C^\cdot({\mathcal U},{\mathcal E}^{1,\cdot}))_{tot})$. Then, we use Lemma \ref{atiyah}.\\
Finally note that $H^1((C^\cdot(\mathcal{U},\mathcal{E}^{1,\cdot}_X)_{tot})
\simeq H^1(\mathcal{U},\Omega^1_X)$ 
and we have $H^1(\mathcal{U},\Omega^1_X)\simeq H^1(X,\Omega^1_X)$ if the covering $\mathcal U$ is 
by Stein open sets.\\
(3) This follows from (2) and Corollary \ref{existence}, or directly: 
Assume that the cohomology class of $\kappa$ in $H^1(H^0(X,{\mathcal E}^{1,*}_X))$ vanishes. 
Then there is $\beta\in H^0(X,{\mathcal E}^{1,0})$ such that 
$\bar{\partial}\alpha_i=\bar{\partial}\beta$ on $U_i$. Then $\alpha_i':=\alpha_i-\beta$ is 
holomorphic on $U_i$, $\alpha_j'-\alpha_i'=\frac{dg_{ij}}{g_{ij}}$, so it defines a complex analytic 
connection on 
$\mathcal{L}$.

\vskip.1in

Altogether, the curvature and the $\bar{\partial}$-curvature decide respectively about whether the given connection (on a differentiable line bundle $\mathcal{E}_X\otimes_{\mathcal{O}_X}\mathcal{L}$) is integrable or complex analytic. Their cohomology classes $-2\pi ic_1(\mathcal{L})_\mathbb{C}$ and $b(\mathcal{L})$ decide respectively about the existence of such connections.  In several cases we have that the vanishing of the complex first Chern class implies the existence of an integrable complex-analytic connection: if $X$ is Stein (Lemma \ref{Stein}) or compact K\"ahler (Lemma \ref{integrable} or Theorem \ref{kahler})
or $X,\mathcal{L}$ algebraic (Theorem \ref{tors}). In these cases we must necessarily have that $b(\mathcal{L})=0$ - an implication which does not hold in general, see Example 4.5. 

If the Atiyah 
obstruction vanishes, according to Lemma \ref{b(L)} the invertible sheaf $\mathcal{L}$ has a connection $\nabla$, but it  is not clear that this connection is an integrable (complex-analytic) connection even if the complex Chern class vanishes 
(see
\S \ref{integ} below).

If we fix a Hermitian metric on $\mathcal{E}_X\otimes_{\mathcal{O}_X}\mathcal{L}$ there is a canonical choice of a connection which respects the complex structure: the metric connection. In this case the $\bar{\partial}$-curvature coincides with the (usual) curvature of the connection, i.e. $d\alpha_i=\bar{\partial}\alpha_i$ above (see \cite{GH} \S  0.5, p. 73). Therefore:

\begin{lemma}\label{curvmet} : 
\begin{enumerate}
\item{} The curvature of a metric connection is a $(1,1)$-form (see \cite{GH} p. 75).
\item{}  In the de Rham cohomology, the complex Chern class 
$c_1(\mathcal{L})$ is representable by a $(1,1)$-form.
\end{enumerate}
\end{lemma}

\subsection{}{\bf $Pic^{an}_c(X)$ and $Pic^{an}_{ci}(X)$} \label{integ}\\

Recall that a connection $\nabla$ is integrable if its curvature vanishes.

When $\mathcal{L}=\mathcal{O}_X$ and $\nabla(f)=df+f\omega$, we have the value of the curvature $R_\nabla$
of the connection $\nabla$:
$$R_{\nabla}=d\omega+ \omega\wedge\omega$$
(see I 3.2.2 of \cite{[D]}, p. 23). In the case of invertible sheaves we have 
$\omega\wedge\omega=0$. 

More generally, recall that a connection is given by a ${\mathbb C}$-linear morphism:
$$\nabla^1:{\mathcal L}\to \Omega^1_X\otimes {\mathcal L}=\Omega^1_X({\mathcal L})$$
It defines a ${\mathbb C}$-linear morphism:
$$\nabla^2:\Omega^1_X({\mathcal L})\to \Omega^2_X({\mathcal L})$$ 
by the formula:
$\nabla^2(\omega\otimes s)=d\omega \otimes s -\omega\wedge \nabla(s)$ (see I (2.4) and (2.9)
of \cite{[D]}).

\begin{definition}
The connection $\nabla=\nabla^1$ is said to be integrable if $\nabla^2\circ\nabla^1=0$.
\end{definition}

In particular the trivial connection $d$ on ${\mathcal O}_X$ is integrable.
As we did for the group $Pic^{an}_c(X)$, the isomorphism classes 
 of analytic invertible sheaves with integrable connection are a group $Pic^{an}_{ci}(X)$ in 
which the neutral element is the class of $({\mathcal O}_X,d)$ and the class of (${\mathcal L}_1,\nabla_1$)
and the one of (${\mathcal L}_2,\nabla_2$) give the class of 
(${\mathcal L}_1\otimes {\mathcal L}_2,\nabla$), where:
$$\nabla(s_1\otimes s_2)=\nabla_1(s_1)\otimes s_2+s_1\otimes \nabla_2(s_2).$$

One can prove (see \cite{[D]} using  Th\'eor\`eme 2.17 Chap. I p. 12) that, if (${\mathcal L}_1,\nabla_1$)
and (${\mathcal L}_2,\nabla_2$) are integrable connections, the connection:
$$({\mathcal L}_1\otimes {\mathcal L}_2,\nabla)$$ 
is also integrable. One can see this directly, too, using that the sum of closed forms is closed.

The curvature of a connection (${\mathcal L},\nabla$) defines an ${\mathcal O}_X$-homomorphism:
$${\mathcal L}
\rightarrow \Omega_X^2\otimes {\mathcal L}$$
which must be of the form
$$s\mapsto \omega\otimes s.$$
Therefore it is given by an element $\omega$ of $H^0(X,\Omega_X^2)$. If this cohomology group vanishes, we have
$Pic^{an}_{ci}(X)\simeq Pic^{an}_{c}(X)$.

\begin{proposition}\label{inj}
Let $X$ be a complex manifold. We have an exact sequence
$$0\rightarrow Pic^{an}_{ci}(X)\rightarrow Pic^{an}_{c}(X)\rightarrow H^0(X,\Omega^2_X).$$
\end{proposition}

\noindent{\bf Proof}. Let $({\mathcal L},\nabla)$ be an integrable connection.

Assume this connection is isomorphic to the trivial connection $({\mathcal O}_X,d)$, the class
of the connection $({\mathcal L},\nabla)$ is therefore the class of the trivial connection. This means
that the map $Pic^{an}_{ci}(X)\rightarrow Pic^{an}_{c}(X)$ is an injection.

The mapping $Pic_c^{an}(X)\to H^0(X,\Omega^2_X)$ associates the curvature of $\nabla$ 
with the isomorphism class of $(\mathcal{L},\nabla)$. It is well-defined: 
if $(\mathcal{L},\nabla)$ and $(\mathcal{L}',\nabla')$ are isomorphic the two connections are represented by the same differential forms with respect to corresponding sections. The exactness at $Pic_c^{an}(X)$ is obvious.
\vskip.1in

In fact the following proposition shows that $Pic^{an}_{ci}(X)$ is of topological nature:

\begin{proposition}\label{iso}
We have the isomorphism:
$$Pic^{an}_{ci}(X)\simeq H^1(X,{\mathbb C}^*).$$
\end{proposition}

\noindent {\bf Proof}. According to Th\'eor\`eme 2.17 in chapter I of \cite{[D]} there is an equivalence 
of categories between the category of line bundles coming from a complex local system together with their canonical
connection with the category of line bundles with an integrable connection.

When $X$ is connected, since the  group ${\mathbb C}^*$ is abelian, the isomorphism classes of invertible sheaves coming from a complex local system together with their canonical
connection is $Hom(\pi_1(X,*),{\mathbb C}^*) $ (using e.g. Corollaire 1.4 of Chap. I of \cite{[D]}).

This gives:
$$Hom(\pi_1(X,*),{\mathbb C}^*) \simeq Pic^{an}_{ci}(X)$$
where $*$ is some point of $X$.
In fact, if we want to associate to a local system a homomorphism $h :\pi_1(X,*)\to {\mathbb C}^*$, 
we have to identify the fiber over $*$ with ${\mathbb C}$, so $h$ is 
unique up to replacement by the conjugate $\alpha\mapsto c\,h(\alpha)\,c^{-1},c\in {\mathbb C}^*$. But
$c\,h(\alpha)\,c^{-1}=h(\alpha)$ because the group ${\mathbb C}^*$ is abelian.

The homology group $H_1(X,{\mathbb Z})$ being the abelianization of the fundamental group
$\pi_1(X,*)$, we have:
$$Hom(H_1(X,{\mathbb Z}),{\mathbb C}^*) \simeq Hom(\pi_1(X,*),{\mathbb C}^*).$$

Furthermore: 
$${\rm Ext}^1(H_0(X,{\mathbb Z}),{\mathbb C}^*)=0,$$
because the abelian group $H_0(X,{\mathbb Z})$ is free, and the Universal coefficient formula implies 
$$H^1(X,{\mathbb C}^*)\simeq Hom(H_1(X,{\mathbb Z}),{\mathbb C}^*).$$

This yields:
$$H^1(X,{\mathbb C}^*)\simeq Pic^{an}_{ci}(X).$$

\vskip.1in
In general, when $X$ has the connected components $X_i,i\in I$, we conclude that:
$$Pic^{an}_{ci}(X)\simeq\prod_{i\in I}Pic^{an}_{ci}(X_i)\simeq\prod_{i\in I}H^1(X_i,{\mathbb C}^*)\simeq H^1(X,{\mathbb C}^*)$$

\vskip.1in
\noindent {\bf Remark.} Alternative proof: We can observe that the group $H^1(X,{\mathbb C}^*)$ 
classifies the local systems 
of one dimensional complex vector spaces on $X$ (see Theorem 3.3 of  \cite{St}), up to isomorphism, 
because the local transition functions are locally constant. The same is true for $Pic^{an}_{ci}(X)$ 
as mentioned above.

\vskip.1in
In particular the isomorphism:
$$Hom(H_1(X,{\mathbb Z}),{\mathbb C}^*) \simeq H^1(X,{\mathbb C}^*)\simeq Pic^{an}_{ci}(X)$$
implies:

\begin{corollary} \label{isom}
Let $f:X\to Y$ be a holomorphic map between two complex manifolds such that it 
induces an isomorphism $H_1(X,{\mathbb Z})\to H_1(Y,{\mathbb Z})$, then:
$$Pic^{an}_{ci}(X)\simeq Pic^{an}_{ci}(Y).$$
\end{corollary}

\subsection{Relation to Deligne cohomology}\label{Deligne} The preceding subsection 
is closely related to special cases of Deligne cohomology. Let us start with the description 
of $Pic^{an}_c\,X$ as a hypercohomology group.

Let $\mathcal{U}$ be an open covering of $X$. Let $Pic^{an}\,{\mathcal U}$ be the group of isomorphism classes of invertible $\mathcal{O}_X$-modules which are trivial on the $U_i$. Let $Pic^{an}_c\,{\mathcal U}$ be the group of isomorphism classes of such sheaves with connection. First, $Pic^{an}\,{\mathcal U}\simeq H^1({\mathcal U},{\mathcal O}^*_X)$. Let 
${\mathcal S}^\cdot$ be the non-negative complex: 
$${\mathcal O}^*_X\stackrel{g\mapsto\frac{dg}{g}}{\longrightarrow}\Omega^1_X\to 0\to\ldots.$$ 
Let $(C^\cdot({\mathcal U},{\mathcal S}^\cdot))_{tot}$ be the total (or the simple) complex 
associated to the bi-graded complex $C^\cdot({\mathcal U},{\mathcal S}^\cdot)$ (see e.g. \cite{Br} p. 14), and 
$\mathbb{H}^1({\mathcal U},{\mathcal S}^\cdot):=\mathbb{H}^1(C^\cdot({\mathcal U},{\mathcal S}^\cdot)_{tot})$
Then:

\begin{lemma} \label{hyp} a) $Pic^{an}_c\,{\mathcal U}\simeq \mathbb{H}^1({\mathcal U},{\mathcal S}^\cdot)$.\\
b) $Pic^{an}_c\,X\simeq \check{\mathbb{H}}^1(X,{\mathcal S}^\cdot)\simeq \mathbb{H}^1(X,{\mathcal S}^\cdot)$ {\rm (cf. \cite{Br} Theorem 2.2.20, p. 80)}.
\end{lemma}

\noindent {\bf Proof:} a) Use Lemma \ref{Cech} (See \ref{remarkp11}).\\
b) Take the direct limit with respect to open coverings $\mathcal U$. The second isomorphism 
holds because $X$ is paracompact (see \cite{Br} Theorem 1.3.13, p. 32).

\vskip.1in
As a consequence, we obtain the exact sequence of Theorem \ref{SE} again:\\
We have an exact sequence of complexes:
$$0\to C^{\cdot+1}({\mathcal U},\Omega^1_X)\to (C^\cdot({\mathcal U},{\mathcal S}^\cdot))_{tot}\to C^\cdot({\mathcal U},{\mathcal O}^*_X)\to 0$$ 
where 
$(C^\cdot({\mathcal U},{\mathcal S}^\cdot))_{tot}$ is the total (or the simple) complex 
associated to the bi-graded complex $C^\cdot({\mathcal U},{\mathcal S}^\cdot)$ (see e.g. \cite{Br} p. 14). Note that $H^1(V,\Omega^1_X)=0$ for $V=U_{i_0}\cap\ldots\cap U_{i_q}$ because $V$ is Stein.\\ 
This exact sequence induces a long exact cohomology sequence
$$\ldots\to H^k({\mathcal U},{\mathcal O}^*_X)\to H^k({\mathcal U},\Omega^1_X)\to \mathbb{H}^{k+1}({\mathcal U},{\mathcal S}^\cdot)\to H^{k+1}({\mathcal U},{\mathcal O}^*_X)\to\ldots$$
After this take the direct limit and replace \v{C}ech (hyper)cohomology by the usual one.

\vskip.1in
Now let us turn to Deligne cohomology. Let us recall its definition (see \cite{EV} p. 45). Put 
$\mathbb{Z}(p):=(2\pi i)^p\mathbb{Z}\subset\mathbb{C}$. Let $\mathbb{Z}(p)_{\mathcal D}$ be the following non-negative complex:
$$\mathbb{Z}(p)_X\to \Omega^0_X\to\ldots\to \Omega^{p-1}_X\to 0\to\ldots$$
where the first arrow is the inclusion. Then the Deligne cohomology $H_{\mathcal D}^*(X,\mathbb{Z}(p))$ is defined as the hypercohomology $\mathbb{H}^*(X,\mathbb{Z}(p)_{\mathcal D})$. 

Looking at the commutative diagram
$$\begin{array}{cccclcl}
\mathbb{Z}(p)_X&\to&\mathcal{O}_X&\to& \Omega^1_X&\to\ldots\to& \Omega^{p-1}_X\\
\downarrow&&\downarrow&&\downarrow\cdot (2\pi i)^{-p+1}&&\downarrow\cdot (2\pi i)^{-p+1}\\
0&\to&{\mathcal O}^*_X&\stackrel{f\mapsto\frac{df}{f}}{\to}&\Omega^1_X&\to\ldots\to& \Omega^{p-1}_X
\end{array}$$
where the second verical arrow is given by $f\mapsto exp((2\pi i)^{-p+1}f)$
we see that the complex above is quasi-isomorphic to
$$0\to {\mathcal O}^*_X\stackrel{f\mapsto\frac{df}{f}}{\to}\Omega^1_X\to\ldots\to \Omega^{p-1}_X\to 0\to\ldots$$
For $p=1$, we obtain that ${\mathbb Z}(1)_{\mathcal D}$ is quasi-isomorphic to 
$ {\mathcal O}^*_X(-1)$, cf. \cite{B} p. 2038, so $H^1_{\mathcal D}(X,\mathbb{Z}(1))\simeq H^0(X,{\mathcal O}^*_X)$ and
$H^2_{\mathcal D}(X,\mathbb{Z}(1))\simeq Pic^{an}(X)$.\\
For $p=2$, we get that ${\mathbb Z}(2)_{\mathcal D}$ is quasi-isomorphic to 
${\mathcal S}^\cdot(-1)$, cf. \cite{EV} p. 46, so $Pic_c^{an}(X)\simeq H^2_{\mathcal D}(X,\mathbb{Z}(2))$ because of Lemma 
\ref{hyp} (see the remark of Deligne quoted in \cite{B} below p. 2039).\\
For $p\ge \dim\,X+1$ the complex is quasi-isomorphic to $0\to{\mathcal O}^*_X\to d{\mathcal O}_X\to 0\to\ldots$, see beginning of subsection \ref{GR}; by Poincar\'e Lemma, it is also quasi-isomorphic to $0\to \mathbb{C}^*_X\to 0\to\ldots$.\\
So $H^2_{\mathcal D}(X,\mathbb{Z}(p))\simeq H^1(X,{\mathbb C}^*_X)\simeq Pic_{ci}^{an}(X)$, using Proposition \ref{iso}.\\
For $p>2$, $H^2_{\mathcal D}(X,\mathbb{Z}(p))$ does not depend on $p$.\\

We obtain altogether, cf. \cite{Ga} p. 156:

\begin{lemma}  a) $H^2_{\mathcal D}(X,\mathbb{Z}(1))\simeq Pic^{an}(X)$.\\
b) $H^2_{\mathcal D}(X,\mathbb{Z}(2))\simeq Pic_c^{an}(X)$.\\
c) $H^2_{\mathcal D}(X,\mathbb{Z}(p))\simeq Pic_{ci}^{an}(X)$ for $p>2$.
\end{lemma}

\subsection{}{$Pic^{an}(X)$ and $Pic^{an}_{ci}(X)$}
\vskip.1in
The first exact sequence of \S \ref{GR} gives a long exact sequence which fits into a commutative diagram:

\begin{theorem}\label{CI}
We have a commutative diagram with exact rows:
$$\footnotesize\begin{array}{ccccccccccccc}0&\to &H^0(X,{\mathbb C}^*_X)&\to &H^0(X,{\mathcal O}^*_X)
&\to&H^0(X,d{\mathcal O}_X)
&\to&Pic^{an}_{ci}(X)&\to &Pic^{an}(X)&\to &H^1(X,d{\mathcal O}_X)\\
&&\downarrow&&\downarrow&&\downarrow&&\downarrow&&\downarrow&&\downarrow\\
0&\to &H^0(X,{\mathbb C}^*_X)&\to &H^0(X,{\mathcal O}^*_X)
&\to&H^0(X,\Omega^1_X)
&\to&Pic^{an}_{c}(X)&\to &Pic^{an}(X)&\to &H^1(X,\Omega^1_X)\\
\end{array}
$$
\end{theorem}

\noindent{\bf Proof.} The exactness of the upper line is consequence of Proposition \ref{iso}. The relevant part of the exactness of the lower line was given in Theorem \ref{SE}, since the vertical map: 
$$H^0(X,d{\mathcal O}_X)\to H^0(X,\Omega^1_X)$$ 
is injective.

\vskip.1in
\noindent {\bf Remark:} We may also argue using hypercohomology: 

In the upper row compare $0\to {\mathcal O}_X^*\to d{\mathcal O}_X\to 0$ with $0\to {\mathcal O}_X^*\to 0$, in the lower row 
$0\to {\mathcal O}_X^*\to \Omega^1_X\to 0$ with $0\to {\mathcal O}_X^*\to 0$.
\vskip.1in
In particular, we observe that:

\begin{lemma}
If the complex manifold $X$ is compact with an invertible $\mathcal{O}_X$-module ${\mathcal L}$ on $X$ and if $\nabla_1$ and $\nabla_2$ 
are two connections on ${\mathcal L}$ such that $({\mathcal L},\nabla_1)\simeq({\mathcal L},\nabla_2)$, we must have $\nabla_1=\nabla_2$.
\end{lemma}

\noindent{\bf Proof.} We have $(\nabla_1-\nabla_2)(s)=\omega\otimes s$ where $\omega\in H^0(X,\Omega^1_X)$ is mapped to $0\in Pic^{an}_c(X)$. So there is $g\in H^0(X,\mathcal{O}^*_X)$ such that $\omega=\frac{dg}{g}$. Since 
$H^0(X,{\mathbb C}^*)=H^0(X,{\mathcal O}^*_X)$ because global functions on $X$ are 
locally constant on a compact space, we have that $\omega=0$.

\begin{lemma}\label{chern}
Let $X$ be a complex manifold (eventually non-compact). An element $x\in H^2(X,{\mathbb Z})$ is sent on $0$ in $H^2(X,{\mathbb C})$ if and 
only if it is the Chern
class of an invertible $\mathcal{O}_X$-module which can be endowed with an integrable connection.
\end{lemma}

\noindent{\bf Proof.} We have a commutative diagram:

$$\begin{array}{ccccccccc}0&\to &{\mathbb Z}&\to&{\mathbb C}&\to&{\mathbb C}^*&\to &0\\
&&\downarrow&&\downarrow&&\downarrow&&\\
0&\to &{\mathbb Z}&\to&{\mathcal O}_X&\to&{\mathcal O}^*_X&\to &0\\
\end{array}$$
with exact rows. This leads to a commutative diagram:
$$\begin{array}{cccc}&H^1(X,{\mathbb C}^*_X)&\to&H^2(X,{\mathbb Z})\\
&\downarrow&&\downarrow\\
&H^1(X,{\mathcal O}^*_X)&\to&H^2(X,{\mathbb Z})\\
\end{array}$$
The lower arrow associates to each invertible sheaf its first Chern class, therefore the upper
arrow associates to each invertible sheaf with an integrable connection the first Chern class of the invertible sheaf. Now consider the upper row of the first diagram. It leads to an exact sequence:
$$H^1(X,{\mathbb C}^*)\to H^2(X,{\mathbb Z})\to H^2(X,{\mathbb C}),$$
which gives our result.\\

\noindent {\bf Remark:} We can make Proposition \ref{inj} more precise: There is an exact sequence
$$0\to Pic_{ci}^{an}(X)\to Pic_c^{an}(X)\to H^0(X,d\Omega^1_X)\to H^2(X,\mathbb{C}^*_X)$$

Compare the non-negative complexes ${\mathcal O}_X^*\to d{\mathcal O}_X\to 0$ and ${\mathcal O}_X^*\to \Omega^1_X\to 0$, see subsection \ref{Deligne}. The cokernel is quasi-isomorphic to $0\to\Omega^1_X/d{\mathcal O}_X\to 0$, i.e. to $0\to d\Omega^1_X\to 0$.\\

\subsection{Compact K\"ahler manifolds} In the case $X$ is a compact K\"ahler manifold, we can apply Hodge Theory.\\
In particular, we can view $H^q(X,\Omega^p_X)$ as a subspace of $H^{p+q}(X;\mathbb{C})$:\\
Let $\mathcal{H}^{p,q}(X)$ be the space of harmonic $(p,q)$-forms, $\mathcal{H}^r(X)$ the one of harmonic $r$-forms.\\
Then $\mathcal{H}^r(X)=\oplus_{p+q=r}\mathcal{H}^{p,q}(X)$.\\
Furthermore, $\mathcal{H}^{p,q}(X)\simeq H^q(X,\Omega^p_X)$ (Dolbeault) and\\
$\mathcal{H}^r(X)\simeq H^r(X;\mathbb{C})$ (de Rham).\\
Let $H^{p,q}(X)$ be the image of $\mathcal{H}^{p,q}(X)$ in $H^{p+q}(X;\mathbb{C})$.\\
Then the first part of the following Lemma is well-known:

\begin{lemma}\label{1,1}
Let $X$ be a compact K\"ahler manifold, $\mathcal{L}$ an invertible sheaf on $X$.\\
a) {\rm (see \cite{GH} Ch. 3.3, p. 417)} The complex Chern class $c_1(\mathcal{L})_\mathbb{C}$ of $\mathcal{L}$ is in $H^{1,1}(X)$.\\
b) {\rm (see \cite{A} Prop. 12, p. 196)} With the identifications above, $b(\mathcal{L})=-2\pi ic_1(\mathcal{L})_\mathbb{C}$.\\
\end{lemma}

\noindent{\bf Proof}. a) As shown in \cite{GH} Ch. 0.7, p. 116, the image of $\mathcal{H}^{p,q}(X)$ in $H^{p+q}_{DR}(X)$ consists of the classes of closed $(p,q)$-forms. By Lemma \ref{curvmet}, the complex first Chern class - viewed as an element of $H^2_{DR}(X)$ - can be represented by a form of type $(1,1)$.  \\
b) We have a commutative diagram with exact rows
$$\begin{array}{cclclccccc}
0&\to&\mathbb{Z}_X&\to&\mathcal{O}_X&\stackrel{f\mapsto e^{2\pi i f}}{\to}&\mathcal{O}^*_X&\to&0\\
&&\downarrow\cdot 2\pi i&&\downarrow\cdot 2\pi i&&\downarrow&&\\
0&\to&\mathbb{C}_X&\to&\mathcal{O}_X&\to&d\mathcal{O}_X&\to&0
\end{array}$$
We get a commutative diagram
$$\begin{array}{ccc}
H^1(X,\mathcal{O}^*_X)&\to& H^2(X;\mathbb{Z})\\
\downarrow&&\downarrow \cdot 2\pi i\\
H^1(X,d\mathcal{O}_X)&\to& H^2(X;\mathbb{C})\\
\downarrow&&\\
H^1(X,\Omega^1_X)&&
\end{array}$$
Look at the images of $(g_{ij})$.\\
By \cite{Hi} Theorem 4.3.1, p. 62, we have that the image in $H^2(X;\mathbb{C})$ is $2\pi i c(\mathcal{L})_\mathbb{C}$.\\
By Lemma \ref{atiyah}, the image in $H^1(X,\Omega^1_X)$ is $-b(\mathcal{L})$.\\
Now $H^1(X,d\mathcal{O}_X)\simeq \mathbb{H}^1(X,\Omega^{\ge 1}_X)\simeq F^1H^2(X;\mathbb{C})$, where $F^p$ denotes the Hodge filtration. Since the first Chern class is real we have that $(\frac{1}{2\pi i}\frac{dg_{ij}}{g_{ij}})$ represents an element in $(F^1\cap\bar{F}^1)(H^2(X;\mathbb{C})$. This group corresponds to $H^1(X,\Omega^1_X)$, so we obtain our statement.\\
Note that the proof of b) in \cite{A} loc. cit. works only if $\dim\,X=1$ because it uses an exact sequence of the form
$$0\to\mathbb{C}_X\to\mathcal{O}_X\to\Omega^1_X\to 0$$

\vskip.1in

\noindent Now we may sharpen Lemma \ref{chern}:

\begin{lemma} \label{integrable} Let $X$ be a compact K\"ahler manifold, $\mathcal L$ an invertible sheaf on $X$. Then the following conditions are equivalent:\\
a) $\mathcal L$ admits an integrable connection, \\
b) $\mathcal L$ admits a connection, \\ 
c) the Chern class of $\mathcal{L}$ is a torsion element.
\end{lemma}

\noindent{\bf Proof.} That the Chern class is a torsion element means that the complex Chern class vanishes, because it is known that the cohomology group $H^2(X,{\mathbb Z})$ is finitely generated when $X$ is compact. \\
a) $\Leftrightarrow$ c): $\mathcal L$ admits an integrable connection if and only if the image of $\mathcal L$ in $H^1(X,d{\mathcal O}_X)$ vanishes, by Theorem \ref{CI}.\\
The composition $Pic^{an}\,X\to H^1(X,d\mathcal{O}_X)\to H^2(X;\mathbb{C})$ is given by $[\mathcal{L}]\mapsto 2\pi i c_1(\mathcal{L})_\mathbb{C}$, see proof of the preceding lemma.\\
To prove a) $\Leftrightarrow$ c) it is therefore sufficient to show that the mapping $H^1(X,d{\mathcal O}_X)\to H^2(X;\mathbb{C})$ is injective. 

The sheaf $d{\mathcal O}_X$ is quasi-isomorphic to the complex:
$$0\to \Omega^1_X\to \Omega^2_X\to \ldots$$ 
translated by $-1$ which is the filtration b\^ete (or Hodge filtration) in degree one of
the De Rham complex. This gives an isomorphism:
$$H^1(X,d{\mathcal O}_X)\simeq F^1H^2(X,{\mathbb C})$$
where $F^1$ gives the first term of the Hodge filtration of $H^2(X,{\mathbb C})$. This implies 
the desired injectivity. \\
Now b) $\Leftrightarrow$ c), because we know that b) holds if and only if $b(\mathcal{L})=0$
by Lemma \ref{b(L)}. The rest follows from the preceding lemma \ref{1,1}.\\

\noindent{\bf Remark.} The proof shows that the complex first
Chern class is contained in $F^1H^2(X,{\mathbb C})$. Since it is real, it must be in $H^{1,1}$. It gives 
another proof of Lemma \ref{1,1} a).\\

In the preceding Lemma \ref{integrable}, in the case of compact K\"ahler manifolds, 
we can sharpen the fact that a) $\Leftrightarrow$ b):

\begin{theorem}\label{kahler}
If $X$ is a compact K\"ahler manifold, a connection on an invertible sheaf is integrable.
\end{theorem}

\noindent{\bf Proof}. From Theorem \ref{CI} we have the following commutative diagram:

$$\begin{array}{ccccccc}
H^0(X,d{\mathcal O}_X)&\to&Pic^{an}_{ci}(X)&\to&Pic^{an}(X)&\to &
I\cr
\downarrow&&\downarrow&&\|&&\downarrow\cr
H^0(X,\Omega^1_X)&\to &Pic^{an}_c(X)&\to&Pic^{an}(X)&\to&H^1(X,\Omega^1_X)\cr
\end{array}$$
where $I$ is the image of the map $(Pic^{an}(X)\to H^1(X,d{\mathcal O}_X))$ and
the upper 
and lower lines are exact.

Since on a compact K\"ahler manifold, global forms are closed, we have 
$H^0(X,d{\mathcal O}_X)=H^0(X,\Omega^1_X)$ (see \cite{GH} second statement
on the top of page 110 and reading ``... is never exact if $\neq 0$")

From the proof of the preceding lemma, we have: $H^1(X,d{\mathcal O}_X)=F^1H^2(X,{\mathbb C})$. 
Therefore, we have $H^1(X,d{\mathcal O}_X)=H^{2,0}(X)\oplus H^{1,1}(X)$. Lemma \ref{1,1} shows that 
$I\subset H^{1,1}(X)$, so the last vertical arrow is injective.

The Five Lemma shows that the map: $Pic^{an}_{ci}(X)\to Pic^{an}_c(X)$ is an epimorphism. We 
already know that it is injective, see Proposition \ref{inj}, therefore, if $X$ is a compact K\"ahler manifold, 
this map is an isomorphism.

Now, consider an invertible sheaf ${\mathcal L}$ with connection $\nabla$. Since 
$Pic^{an}_{ci}(X)\simeq Pic^{an}_c(X)$, there is an invertible sheaf ${\mathcal L}_1$ with an  integrable connection $\nabla_1$ such that $({\mathcal L},\nabla)$ is isomorphic to $({\mathcal L}_1,\nabla_1)$. But this means that $\nabla$ is integrable, too. This proves our Theorem.

Alternative proof: Let $\nabla$ be a connection on $\mathcal{L}$. By Lemma \ref{1,1}, we have $c_1(\mathcal{L})_\mathbb{C}=0$, because $b(\mathcal{L})=0$. The curvature of $\nabla$ is an element of $H^0(X,\Omega^2_X)$. Since $H^0(X,\Omega^2_X)\subset H^2(X;\mathbb{C})$ we obtain that the curvature is $0$. So the connection $\nabla$ is integrable.\\

\section{Algebraic case} 

\subsection{} Suppose now that $X$ is a smooth complex algebraic variety. The underlying analytic space 
$X^{an}$ is a paracompact complex manifold. One has an analogue of Theorem \ref{CI} but one has to be careful with 
the upper row because one has no longer a Poincar\'e lemma. In fact we have to replace the sheaf 
$d{\mathcal O}_X$ by the sheaf $^c\Omega^1_X$ of closed Pfaffian forms on $X$.\\

\begin{theorem} \label{alg} Let $X$ be a smooth complex algebraic variety. Then we have a commutative diagram with exact rows
\footnotesize
$$\begin{array}{ccccccccccccc}
0&\to&H^0(X,\mathbb{C}^*_X)&\to&H^0(X,{\mathcal O}^*_X)&\to&H^0(X,{^c\Omega}^1_X)&\to&Pic_{ci}(X)&\to&Pic(X)&\to &H^1(X,{^c\Omega}^1_X)\\
&&\downarrow&&\downarrow&&\downarrow&&\downarrow&&\downarrow&&\downarrow\\
0&\to&H^0(X,\mathbb{C}^*_X)&\to&H^0(X,{\mathcal O}^*_X)&\to&H^0(X,\Omega^1_X)&\to&Pic_c(X)&\to&Pic(X)&\to &H^1(X,\Omega^1_X)
\end{array}$$
\normalsize
\end{theorem}

\noindent{\bf Proof.} We can no longer use the exact sequence of the beginning of section \ref{GR}. Therefore we must proceed in a different way.

Let us check first that the lower row is exact.

Note that the sequence of sheaves: $0\to \mathbb{C}^*_X\to {\mathcal O}^*_X\to \Omega^1_X$ 
is exact. In fact:

Suppose that $h\in{\mathcal O}^*_{X,x}$, where $x$ is a closed point of $X$, $\frac{dh}{h}=0$: Then $h^{an}\in {\mathcal O}^*_{X^{an},x}$ is mapped to $0\in\Omega^1_{X^{an},x}$, so $h^{an}$ is constant, which implies that $h$ is constant.\\
Therefore the sequence: 
$$0\to H^0(X,\mathbb{C}^*_X)\to H^0(X,{\mathcal O}^*_X)\to H^0(X,\Omega^1_X)$$ 
is exact.

The rest goes as in the proof of Theorem \ref{SE}.

The upper row is treated in an analogous way. Note that the connection $\nabla$ on ${\mathcal O}_X$: $$\nabla(f)=df+f\omega$$ 
is integrable if and only if $\omega$ is closed, because the curvature of $\nabla$ is $d\omega$.

\vskip.1in
\noindent Proposition \ref{inj} has an algebraic counterpart:

\begin{proposition}\label{injalg}
Let $X$ be a non-singular complex algebraic variety. We have an exact sequence
$$0\rightarrow Pic_{ci}(X)\rightarrow Pic_{c}(X)\rightarrow H^0(X,\Omega^2_X).$$
\end{proposition}
 
\noindent The proof is similar to the one of Proposition \ref{inj}.

\vskip.1in
\noindent \label{remarkp11} Similarly as in the analytic case (see \S \ref{Deligne}) we can observe that $Pic_c(X)$ is isomorphic to the first \v{C}ech hypercohomology $\check{\H}^1(X,{\mathcal S}^\cdot)$
of the complex ${\mathcal S}^\cdot$:
$$0\to {\mathcal O}^*_X\to \Omega^1_X\to 0 \to \ldots$$
on $X$. 

In fact, let ${\mathcal U}$ be a covering of $X$ by open Zariski subsets of $X$. An element of 
${\H}^1({\mathcal U},{\mathcal S}^\cdot)$ is given by an element 
$((\omega_i),(g_{ij}))\in C^0({\mathcal U},\Omega^1_X)\oplus C^1({\mathcal U},{\mathcal O}^*_X)$
such that $(g_{ij})$ is a cocycle, i.e. $g_{ij}=g_{ik}g_{kj}$ on $U_i\cap U_j\cap U_k$, 
and $\omega_j-\omega_i={dg_{ij}\over g_{ij}}$ on $U_i\cap U_j$.

Assume now that $\mathcal L$ is an invertible $\mathcal{O}_X$-module on $X$ which is endowed 
with a connection $\nabla$. There is a Zariski open covering $\mathcal U$ of $X$ such that for each 
$U_i$ we have a trivialization of ${\mathcal L}|U_i$. Then $\mathcal L$ is represented by some 
cocycle $(g_{ij})$ in $C^1({\mathcal U},{\mathcal O}^*_X)$, and $\nabla|U_i$ corresponds to 
a connection $g\mapsto dg+g\omega_i$ on ${\mathcal O}_{U_i}$. Then 
$\omega_j-\omega_i=\frac{dg_{ij}}{g_{ij}}$ on $U_i\cap U_j$, so we obtain an element of 
${\H}^1({\mathcal U},{\mathcal S}^\cdot)$, hence of $\check{\H}^1(X,{\mathcal S}^\cdot)$.

On the other hand, an element of $\check{\H}^1(X,{\mathcal S}^\cdot)$ comes from an element 
of ${\H}^1({\mathcal U},{\mathcal S}^\cdot)$ which is represented by a cocycle $(g_{ij})$ 
and $(\omega_i)$ for a suitable open Zariski covering ${\mathcal U}$ of $X$. Then $(g_{ij})$ 
defines an invertible $\mathcal{O}_X$-module $\mathcal L$, and 
$(\omega_i)$ defines a connection on $\mathcal L$.

Therefore $Pic_c(X)\simeq \check{\mathbb{H}}^1(X,{\mathcal S}^\cdot)$.

The proof of the preceding Theorem gives the following exact sequence:
$$\check{H}^0(X,{\mathcal O}^*_X)\to\check{H}^0(X,\Omega^1_X)\to \check{\mathbb{H}}^1(X,{\mathcal S}^\cdot)\to \check{H}^1(X,{\mathcal O}^*_X)\to\check{H}^1(X,\Omega^1_X)$$
Now this sequence can be mapped to the analogous exact sequence with $H$ instead of $\check{H}$ (resp. ${\H}$ instead of $\check{\H}$): look at the long exact hypercohomology sequence of the exact sequence of sheaf complexes:
$$0\to \Omega^1_X\{1\}\to\mathcal{S}^\cdot\to\mathcal{O}^*_X\{0\}\to 0$$
where, for any sheaf $\mathcal T$, the complex $\mathcal T\{k\}$ denotes the complex ${\mathcal T}^\cdot$ with ${\mathcal T}^l=\mathcal T$ for $l=k$ and $=0$ otherwise. 

Now in the case of sheaves we have isomorphisms $\check{H}^k\to H^k$ for $k=0,1$, see \cite{Go} II 5.9 Corollaire, p. 227 (note that $X$ is not paracompact and that we are not only dealing with coherent algebraic sheaves!). By the Five Lemma we obtain that $\check{\mathbb{H}}^1(X,{\mathcal S}^\cdot)\simeq \mathbb{H}^1(X,{\mathcal S}^\cdot)$.
We have proved:
\begin{lemma}
If $X$ is a non-singular complex variety, we have: 
$$Pic_c(X)\simeq \check{\mathbb{H}}^1(X,{\mathcal S}^\cdot)
\simeq \mathbb{H}^1(X,{\mathcal S}^\cdot).$$
\end{lemma}

\vskip.1in
\noindent This motivates a different proof of the exactness of the lower row of Theorem \ref{alg}:

\noindent It is sufficient to prove that 
$\check{\mathbb{H}}^1(X,{\mathcal S}^\cdot)\simeq \mathbb{H}^1(X,{\mathcal S}^\cdot)$. 

\vskip.1in
\noindent In fact, we have:

\noindent \begin{proposition}  Let $X$ be a topological space and $\mathcal{S}^\cdot$ a non-negative complex of abelian groups on $X$. Then the homomorphism
$\check{\mathbb{H}}^k(X,\mathcal{S}^\cdot)\to \mathbb{H}^k(X,\mathcal{S}^\cdot)$ is bijective for $k\le 1$ and injective for $k=2$.
\end{proposition}

\noindent {\bf Proof:} Recall the definition of \v{C}ech hypercohomology, cf. \cite{Br} p. 28: Let $\mathcal{U}$ be an open covering of $X$. Then we may look at the double complex $C^\cdot (\mathcal{U},\mathcal{S}^\cdot)$. The $p$-th hypercohomology of the associated simple complex is denoted by $\mathbb{H}^p(\mathcal{U},\mathcal{S}^\cdot)$. Then $\check{\mathbb{H}}^p(X,\mathcal{S}^\cdot):=\lim\limits_{\to\mathcal{U}}\mathbb{H}^p(\mathcal{U},\mathcal{S}^\cdot)$. Here we take all open coverings $\mathcal{U}$. \\
If we restrict to open coverings $\mathcal{U}=(U_x)_{x\in X}$ such that $x\in U_x$ for all $x$ and use the partial order $\mathcal{U}\le \mathcal{V}$ if and only if $U_x\subset V_x$ for all $x$ 
we can introduce the \v{C}ech double complex 
$\check{C}^\cdot(X,\mathcal{S}^\cdot):=\lim\limits_{\to\mathcal{U}}C^\cdot (\mathcal{U},\mathcal{S}^\cdot)$, similarly as in \cite{Go} p. 223. Then  
$\check{\mathbb{H}}^p(X,\mathcal{S}^\cdot)$ is the $p$-th cohomology group of the associated simple complex.\\
We may use sheaves instead:\\
$\check{C}^\cdot(X,\mathcal{S}^\cdot)=\Gamma(X,\check{\mathcal{C}}^\cdot(\mathcal{S}^\cdot))$. \\
Let $\mathcal{J}^{\cdots}$ be a flabby resolution of $\check{\mathcal{C}}^\cdot(\mathcal{S}^\cdot)$, this is a triple complex, where the third index is the one given by the resolution:
$$0\to \check{\mathcal C}^p({\mathcal S}^q)\to \mathcal{J}^{pq0}\to \mathcal{J}^{pq1}\to \mathcal{J}^{pq2}\to\ldots.$$
 By amalgamation of the first and third index we get a flabby double complex $\mathcal{I}^{\cdot\cdot}$. We have a homomorphism
$\Gamma(X,\check{\mathcal{C}}^\cdot(\mathcal{S}^\cdot))\to\Gamma(X,\mathcal{I}^{\cdot\cdot})$ of double complexes by using:
$$\check{\mathcal C}^p({\mathcal S}^q)\to \mathcal{J}^{pq0}\to\oplus_{r+s=p}{\mathcal J}^{rqs}={\mathcal I}^{pq}$$
For the cohomology of the corresponding total complexes we obtain a homomorphism:
$$\check{\mathbb{H}}^p(X,\mathcal{S}^\cdot)\to \mathbb{H}^p(X,\mathcal{S}^\cdot)$$
because, from Th\'eor\`eme 5.2.1 of \cite{Go}, one has that
$\mathcal{S}^\cdot\to(\check{\mathcal{C}}^\cdot(\mathcal{S}^\cdot))_{tot}$ is a quasi-isomorphism
and
$(\check{\mathcal{C}}^\cdot(\mathcal{S}^\cdot))_{tot}\to(\mathcal{I}^{\cdot\cdot})_{tot}$ is a quasi-isomorphism
by definition. Therefore the following hypercohomologies are isomorphic:
$$\mathbb{H}^p(X,\mathcal{S}^\cdot)\to \mathbb{H}^p(X,(\mathcal{I}^{\cdot\cdot})_{tot})$$
and from the definition of ${\mathcal I}^{\cdot\cdot}$ we have the homomorphism:
$$\check{\mathbb{H}}^p(X,\mathcal{S}^\cdot)\to \mathbb{H}^p(X,(\mathcal{I}^{\cdot\cdot})_{tot})$$
which gives the above homomorphism.
\\
We have a spectral sequence $\check{E}_2^{pq}\Rightarrow \check{\mathbb{H}}^{p+q}(X,\mathcal{S}^\cdot)$:\\ 
Let $\mathcal{U}$ be an open covering of $X$. By \cite{Go} p. 211 the double complex
$\Gamma(X,\mathcal{C}^\cdot (\mathcal{U},\mathcal{S}^\cdot))=C^\cdot (\mathcal{U},\mathcal{S}^\cdot)$ leads to a spectral sequence
$E(\mathcal{U})_2^{pq}=H^p(\mathcal{U},\mathcal{H}^q(\mathcal{S}^\cdot))\Rightarrow \mathbb{H}^{p+q}(\mathcal{U},\mathcal{S}^\cdot)$.\\
By taking the direct limit we arrive at $\check{E}_2^{pq}=\check{H}^p(X,\mathcal{H}^q(\mathcal{S}^\cdot))$.\\
Similarly we have $E_2^{pq}\Rightarrow \mathbb{H}^{p+q}(X,\mathcal{S}^\cdot)$,
where $E_2^{pq}=H^p(X,\mathcal{H}^q(\mathcal{S}^\cdot))$: see \cite{Go} p. 178.\\
Finally we have a homomorphism $\check{E}_r^{pq}\to E_r^{pq}$. If $r=2$ it is bijective for $p+q\le 1$ and injective for $p+q=2$, see \cite{Go} p. 227.\\
By computation we obtain that the same holds for $r=3, 4$. Since $\check{E}_4^{pq}=\check{E}_\infty^{pq}$ and $E_4^{pq}=E_\infty^{pq}$ for $p+q\le 2$ we get the same for $r=\infty$. This implies our proposition.\\

\noindent We can proceed in the same way to prove the exactness of the upper line of the 
diagram of Theorem \ref{alg} by replacing $\Omega^1_X$
by $^c\Omega^1_X$. See Remark after Theorem \ref{CI}.\\

\noindent We have special cases:

\begin{lemma}\label{special1} Let $X$ be complete, $\mathcal L$ an invertible $\mathcal{O}_X$-module on $X$.\\
a) $\mathcal L$ admits an integrable connection if and only if $c_1({\mathcal L})$ is a torsion element.\\
b) Every connection on $\mathcal L$ is integrable.
\end{lemma}

\noindent {\bf Proof:} If $X$ is projective we know that $X^{an}$ is compact K\"ahler, so the result follows by GAGA (see \cite{S} and also \cite{M} p. 152f.) and Lemma \ref{integrable}, Theorem \ref{kahler}. \\
In general we know by \cite{[D3]} \S 5 that we can still apply Hodge theory, so Lemma \ref{integrable} and Theorem \ref{kahler} still hold.\\

For part a) of the lemma it will turn out that the hypothesis that $X$ is complete is unnecessary, see Corollary \ref{tors1} below. For b) we must in general restrict to regular connections, see below (Theorem \ref{int}).\\

Remember that compact K\"ahler manifolds are not automatically algebraic, cf. the case of complex tori, see \cite{Mu} Cor. p. 35.

\begin{lemma}\label{special2} Let $X$ be affine. Then every invertible $\mathcal{O}_X$-module on $X$ admits a connection.
\end{lemma}

\noindent {\bf Proof:} Obvious from Theorem \ref{alg}.

\subsection{Regularity}
It is useful to take the notion of regularity into account. 

\vskip.1in
The regularity has been introduced by P.Deligne in \cite{[D]} Chap II \S 4. By commodity we define 
here the regularity of  integrable connections on an invertible sheaf:

\begin{definition} Let $\mathcal L$ be an invertible $\mathcal{O}_X$-module and $\nabla$ an integrable 
connection on $\mathcal L$. Then $\nabla$ is called regular if there exists a smooth compactification 
$\bar{X}$ of $X$ such that $D:=\bar{X}\setminus X$ is a divisor with normal crossings and that, for all 
$x\in D$, there exists an open Zariski neighbourhood $V$ of $x$ and 
there exists $s\in H^0(V,j_*{\mathcal L}), s$ nowhere vanishing on $V\setminus D$, such that $\nabla s=\alpha\otimes s$ with 
$\alpha\in H^0(V,\Omega^1_{\bar{X}}(log\,D))$. Here $j:X\to\bar{X}$ is the inclusion.
\end{definition}

Note that we can replace:\\ 
``there exists $s\in H^0(V,j_*{\mathcal L}), s$ nowhere vanishing on $V\setminus D$, such that $\nabla s=\alpha\otimes s$"\\
by \\
``for any $s\in H^0(V,j_*{\mathcal L}), s$ nowhere vanishing on $V\setminus D$, we have $\nabla s=\alpha\otimes s$". \\
Here it is important that we deal with invertible sheaves!\\

As P. Deligne noticed, the notion of regularity does not depend on the compactification of $X$
such that the divisor at $\infty$ is a normal crossing divisor (see \cite{[D]} p. 90).

We can define the Picard group $Pic_{cir}X$ of regular integrable connections in an obvious way.

\begin{lemma}\label{exact} There is an exact sequence:
$$
H^0(X,\mathcal{O}^*_X)\to H^0(\bar{X},{^c\Omega}^1_{\bar{X}}(log\,D))\to Pic_{cir}(X)\to Pic(X)\to 
H^1(\bar{X},{^c\Omega}^1_{\bar{X}}(log\,D))$$
\end{lemma}
\noindent {\bf Proof:} The proof is analogous to the proof of Theorem \ref{SE}. 

We first observe that, for any invertible $\mathcal{O}_X$-module ${\mathcal L}$, there is 
a Zariski open covering 
${\mathcal U}=(\bar{U}_i)$ of $\bar{X}$ such that the restriction of ${\mathcal L}$ to 
$U_i=\bar{U}_i\setminus D$
is trivial. 

For this, we may assume that $X$ is connected, hence irreducible. One considers a non-empty and therefore dense Zariski open  subspace $U$ of $X$ on which ${\mathcal L}$ is trivial. 
On $U$, the restriction ${\mathcal L}|U$ has a nowhere vanishing section $s$. This section extends as a
rational section $s_1$ of $\mathcal L$. Let $D_1$ be the divisor of this section - this makes sense because $\mathcal L$ is locally trivial. Now $D_1$ extends to a divisor $\bar{D}_1$ on $\bar{X}$. For any $x\in\bar{X}$ there is an open affine neighbourhood $\bar{V}$ such that $\bar{D}_1|\bar{V}$ is a principal divisor, i.e. divisor of some rational function $\phi_x$. Then $\phi_x^{-1}s_1$ is a nowhere vanishing section of $\mathcal{L}|V$ with $V:=\bar{V}\setminus D$; it gives a trivialization of $\mathcal{L}|V$.

The first arrow is induced by the homomorphism $j_*\mathcal{O}^*_X\to {^c\Omega}^1_{\bar{X}}(log\,D)$ which is defined as follows. Locally, a section $g$ of $j_*\mathcal{O}^*_X$ is of the form $h^{-1}\tilde{g}$, where $h,\tilde{g}$ are regular functions which do not vanish inside $X$. Then the image is defined to be $\frac{dg}{g}=\frac{d\tilde{g}}{\tilde{g}}-\frac{dh}{h}$ which is indeed a closed logarithmic form.

Assume now that $g\in H^0(X,\mathcal{O}^*_X)$ is given. Then the image in $Pic_{cir}(X)$ is given by $\mathcal{O}_X$, together with the connection $f\mapsto df+\frac{dg}{g}$. This is isomorphic to $\mathcal{O}_X$, together with the connection $f\mapsto df$, so we have the trivial element of $Pic_{cir}(X)$.\\
On the other hand, suppose that $\omega\in H^0(\bar{X},{^c\Omega}^1_{\bar{X}}(log\,D))$ is mapped onto the trivial element of $Pic_{cir}(X)$. Then there is a $g\in H^0(X,\mathcal{O}^*_X)$ such that $\omega=\frac{dg}{g}$.\\
This shows the exactness at $H^0(\bar{X},{^c\Omega}^1_{\bar{X}}(log\,D))$.

Then, an element of 
$Pic\;X$ is represented by a cocycle $(g_{ij})$ on a covering ${\mathcal U}$ as defined before. 
This covering comes from an affine covering $\bar{\mathcal U}$ of $\bar{X}$,
where each $g_{ij}$ extends as a rational function with poles inside
$D$ which is a regular and non-vanishing function on $\bar{U}_i\cap \bar{U}_j\setminus D$. Then $\frac{dg_{ij}}{g_{ij}}$ is a closed logarithmic form on $\bar{U}_i\cap\bar{U}_j$: After refining $\mathcal{U}$ if necessary we may assume that we can write 
$g_{ij}=h_{ij}^{-1}\tilde{g}_{ij}$ where $h_{ij}$ and $\tilde{g}_{ij}$ are regular on $\bar{U}_i\cap\bar{U}_j$ and without zeroes in $U_i\cap U_j$.  Then: 
$$\frac{dg_{ij}}{g_{ij}}=\frac{d\tilde{g}_{ij}}{\tilde{g}_{ij}}-\frac{dh_{ij}}{h_{ij}}$$ 
is a closed logarithmic form. This defines the map: 
$$Pic\,X\to H^1(\bar{X},{^c\Omega}^1_{\bar{X}}(log\,D)).$$
On the other hand, a regular integrable connection 
on ${\mathcal O}_X$ is of the form $g\mapsto dg+\omega g$ with
$\omega\in H^0(\bar{X},^c\Omega^1_{\bar{X}}(log\,D))$, i.e. the map from
$H^0(\bar{X},^c\Omega^1_{\bar{X}}(log\,D))$ into $Pic_{cir}(X)$ is given by:
$$\omega\mapsto ({\mathcal O}_X,\nabla)$$
where $\nabla(g)=dg + \omega g$ . Then, the composition:
$$H^0(\bar{X},^c\Omega^1_{\bar{X}}(log\,D))\to Pic_{cir}(X)\to Pic\,X$$
is zero. Let $({\mathcal L},\nabla)$ a regular integrable connection on the invertible $\mathcal{O}_X$-module ${\mathcal L}$ where ${\mathcal L}$
is isomorphic to ${\mathcal O}_X$. The pair $({\mathcal L},\nabla)$ is isomorphic to $({\mathcal O}_X,\nabla_0)$ for some connection $\nabla_0$, and there is a closed logarithmic form 
$\omega\in H^0(\bar{X},^c\Omega^1_{\bar{X}}(log\,D))$, such that 
$\nabla_0(g)=
dg+\omega g$. This proves the exactness of the sequence at $Pic_{cir}(X)$.

Now let an element of $Pic\,X$ whose image in $H^1(\bar{X},{^c\Omega}^1_{\bar{X}}(log\,D))$
is trivial. Such an element is given by an affine covering ${\mathcal U}$ and a
cocycle 
$$({dg_{ij}\over g_{ij}})$$
such that:
$${dg_{ij}\over g_{ij}}= \omega_j-\omega_i$$
where $\omega_i$ is a closed form in ${^c\Omega}^1_{\bar{X}}(log\,D)$ over the Zariski open 
subset $U_i$ of $X$.

As we did in the proof of Theorem \ref{SE}, the element $(\omega_i)$ defines a regular integrable
connection $\nabla$ on an invertible $\mathcal{O}_X$-module ${\mathcal L}$ and the image of the isomorphism class of $({\mathcal L},\nabla)$ is the
element of $H^1(\bar{X},{^c\Omega}^1_{\bar{X}}(log\,D))$ given by the cocycle
$({dg_{ij}\over g_{ij}})$.

It remains to prove that the composition:
$$Pic_{cir}(X)\to Pic\,X\to H^1(\bar{X},{^c\Omega}^1_{\bar{X}}(log\,D))$$
is zero. As in the proof of Theorem \ref{SE}, an element of $Pic_{cir}(X)$ is given by 
$({\mathcal L}|U_i,\nabla |U_i)_i$ such that $({\mathcal L}|U_i,\nabla |U_i)$ is isomorphic 
over the Zariski open subspace $U_i$ to $({\mathcal O}_{U_i},\tilde{\nabla}_i)$ where:
$$\tilde{\nabla}_i(f)=df+\omega_i f$$
for some $\omega_i\in H^0(U_i,{^c\Omega}^1_{\bar{X}}(log\,D))$, 
and, if the element $(g_{ij})$ is the cocycle which defines ${\mathcal L}$, we have:
$${dg_{ij}\over g_{ij}}=\omega_j-\omega_i.$$ 
Since the forms 
$\omega_i$ are closed, reasoning as in the proof of Theorem \ref{SE}, we obtain our assertion.

\vskip.1in
\noindent{\bf Remarks.} 1. In fact, at the beginning we have shown that $j_*{\mathcal L}$ is an invertible $j_*{\mathcal O}_X$-module, $j:X\to\bar{X}$ being the inclusion.\\
2. Again we can prove the lemma by showing that $Pic_{cir}X\simeq \check{\H}^1(X,{\mathcal T}^\cdot)\simeq {\H}^1(X,{\mathcal T}^\cdot)$, where ${\mathcal T}^\cdot$ is the non-negative complex 
$$j_*{\mathcal O}^*_X\stackrel{g\mapsto \frac{dg}{g}}{\longrightarrow} {^c\Omega}^1_{\bar{X}}(log\,D)\longrightarrow 0\longrightarrow\ldots$$ with $j:X\to\bar{X}$ being the inclusion.

\vskip.1in
\noindent In this context it is useful to have:

\begin{lemma} $Pic\,X\simeq H^1(\bar{X},j_*\mathcal{O}^*_X)$.
\end{lemma}

\noindent {\bf Proof:} It is sufficient to show that $R^1j_*\mathcal{O}^*_X=0$. An element of $(R^1j_*\mathcal{O}^*_X)_x$ is represented by an element of $H^1(U\cap X,\mathcal{O}^*_X)$, where $U$ is an open neighbourhood of $x$. After shrinking $U$ if necessary the latter comes from an element of $H^1(U,\mathcal{O}^*_{\bar{X}})$, 
so from a line bundle on $U$. On some smaller neighbourhood $V$ the line bundle is trivial, so the image in 
$H^1(V\cap X,\mathcal{O}^*_X)$ is $0$, which leads to the desired result.

\vskip.1in
\begin{theorem}\label{tors} Let $\mathcal L$ be an invertible $\mathcal{O}_X$-module on $X$. Then $\mathcal L$ admits a regular integrable connection if and only if its first Chern class is a torsion element.
\end{theorem}

\noindent {\bf Proof:} Since the integral cohomology of $X$ is an abelian group
of finite type, the implication $\Rightarrow$ is proved by Lemma \ref{chern}.

Now, consider the implication $\Leftarrow$. 

Suppose that $c_1({\mathcal L})=0$.

Let $\bar{X}$ be a smooth compactification of $X$ such that 
$D:=\bar{X}\setminus X$ is a normal crossing divisor. Suppose that $D$ has $r$ irreducible components. Then 
$\mathcal L$ extends to an algebraic invertible sheaf ${\mathcal L}'$ on $\bar{X}$ with Chern
class $c_1({\mathcal L}')=0$. 

To prove this, we consider the diagram with exact rows:
$$\begin{array}{ccccccc}
\mathbb{Z}^r&\to&Pic\,\bar{X}&\to&Pic\,X&\to&0\\
\downarrow\simeq&&\downarrow c_1&&\downarrow c_1&&\\
H^2(\bar{X}^{an},X^{an};\mathbb{Z})&\to&H^2(\bar{X}^{an};\mathbb{Z})&\to&H^2(X^{an};\mathbb{Z})&&
\end{array}$$
Let $[{\mathcal L}]$ be the class of ${\mathcal L}$. We have assumed that its Chern class is $c_1({\mathcal L})=0$.
Let ${\mathcal L}_1$ be a invertible ${\mathcal O}_{\bar X}$-module whose class has its image equal
to $[{\mathcal L}]$. The Chern class of ${\mathcal L}_1$ comes from an element of 
$H^2(\bar{X}^{an},X^{an};\mathbb{Z})$ which corresponds to an element of $\mathbb{Z}^r$ whose image in 
$Pic\,\bar{X}$ is ${\mathcal L}_2$ which has the same Chern class as ${\mathcal L}_1$. 
The invertible sheaf ${\mathcal L}':={\mathcal L}_1\otimes {\mathcal L}_2^{-1}$ 
has a Chern class equal to $0$ and it extends ${\mathcal L}$.

On the complete non-singular variety $\bar{X}$ we have obtained an invertible sheaf ${\mathcal L}'$ which extends 
${\mathcal L}$ and has Chern class $c_1({\mathcal L}')=0$. By Lemma \ref{special1} the invertible sheaf 
${\mathcal L}'$ is endowed with a integral connection $\nabla'$. The restriction of $\nabla'$ to ${\mathcal L}$
is a regular integral connection.

If $c_1({\mathcal L})=c$, $c$ being a torsion element, by Lemma \ref{chern} 
there is an invertible sheaf with integrable connection on ${X}^{an}$ having $c$ as first Chern class. 
By Deligne's existence theorem (Th\'eor\`eme 5.9 Chap. II of \cite{[D]} p. 97)
we can find an invertible sheaf ${\mathcal L}_1$ on ${X}$ with an integrable connection 
such that $c_1({\mathcal L}_1)=c$. Now 
$c_1({\mathcal L}\otimes({\mathcal L}_1)^{-1})=0$, so by the preceding result
there is an integrable connection on 
${\mathcal L}\otimes({\mathcal L}_1)^{-1}$. So we get an integrable connection on 
${\mathcal L}={\mathcal L}_1\otimes({\mathcal L}'\otimes({\mathcal L}_1)^{-1})$, too.

Therefore if the Chern class of ${\mathcal L}$ is a torsion  element, the invertible sheaf 
$\mathcal L$ yields a regular integrable connection.

\vskip.1in

\begin{corollary}\label{tors1}  Let $\mathcal L$ be an invertible $\mathcal{O}_X$-module. Then the following
conditions are equivalent:
\begin{enumerate}
\item{} $\mathcal L$ admits a regular integrable connection;
\item{} $\mathcal L$ admits an integrable connection;
\item{} ${\mathcal L}^{an}$ admits an analytic integrable connection;
\item{} the Chern class $c_1({\mathcal L})$  of $\mathcal L$ is a torsion element.
\end{enumerate}
\end{corollary} 

\vskip.1in

\subsection {\bf Remark on integrability and regularity.} One may define a notion of regularity for connections which does not suppose 
that the connection is integrable - at least in the case of invertible sheaves.

\begin{definition} Let $\mathcal L$ be an invertible $\mathcal{O}_X$-module and $\nabla$ a
connection on $\mathcal L$. Then $\nabla$ is called regular if there exists a smooth compactification 
$\bar{X}$ of $X$ such that $D:=\bar{X}\setminus X$ is a divisor with normal crossings and that, for all 
$x\in D$ there exists an affine neighbourhood $V$ of $x$ and 
there exists $s\in H^0(V,j_*{\mathcal L})$ which does not vanish on $D$, such that $\nabla s=\alpha\otimes s$ with 
$\alpha\in H^0(V,\Omega^1_{\bar{X}}(log\,D))$. Here $j:X\to\bar{X}$ is the inclusion.
\end{definition}

We can define the group $Pic_{cr}X$  of isomorphism classes of invertible $\mathcal{O}_X$-modules with a regular connection in an obvious way.

In fact, such a regular connection is automatically integrable, because we have:

\begin{theorem}\label{int} If $\mathcal L$ is an invertible $\mathcal{O}_X$-module, every regular connection 
on $\mathcal L$ is integrable. 
\end{theorem}

\noindent {\bf Proof:} We proceed as in the proof of Theorem \ref{kahler}. 

First we show that the mapping: 
$$Pic_{cir}X\to Pic_{cr}X$$ 
is surjective. In fact, we have the following Lemma:

\begin{lemma}\label{log} There is a commutative diagram with exact rows
$$\begin{array}{ccccccccc}
H^0(X,\mathcal{O}^*_X)&\to&H^0(\bar{X},{^c\Omega}^1_{\bar{X}}(log\,D))&\to&Pic_{cir}(X)&\to&Pic(X)&\to&
H^1(\bar{X},{^c\Omega}^1_{\bar{X}}(log\,D))\\
\downarrow&&\downarrow&&\downarrow&&\downarrow&&\downarrow\\
H^0(X,\mathcal{O}^*_X)&\to&H^0(\bar{X},\Omega^1_{\bar{X}}(log\,D))&\to&Pic_{cr}(X)&\to&Pic(X)&\to&
H^1(\bar{X},\Omega^1_{\bar{X}}(log\,D))
\end{array}$$
\end{lemma}

\noindent {\bf Proof:} As in Lemma \ref{exact}, the proof is analogous to the proof of Theorem \ref{SE}. 

The upper line is exact, as we saw in 
Lemma \ref{exact}. Concerning the lower row, we define the map $Pic\,X\to H^1(\bar{X},\Omega^1_{\bar{X}}(log\,D))$ as the composition $Pic\,X\to H^1(\bar{X},{^c\Omega}^1_{\bar{X}}(log\,D))\to H^1(\bar{X},\Omega^1_{\bar{X}}(log\,D))$. 

The map $H^0(\bar{X},\Omega^1_{\bar{X}}(log\,D))\to Pic_{cr}(X)$ is given by:
$$\omega\mapsto ({\mathcal O}_X,\nabla)$$
where the connection $\nabla$ is defined by $\nabla(f)=df +\omega f$. This defines a connection on 
${\mathcal O}_X$ which is regular since $\omega\in H^0(\bar{X},\Omega^1_{\bar{X}}(log\,D))$. 
Therefore, the composition:
$$H^0(\bar{X},\Omega^1_{\bar{X}}(log\,D))\to Pic_{cr}(X)\to Pic\,X$$
is zero. 

Let $({\mathcal L},\nabla)$ be an invertible sheaf with a regular connection whose image is zero in $Pic(X)$. Then 
${\mathcal L}$ is isomorphic to the trivial invertible sheaf ${\mathcal O}_X$ and there is a connection 
$\nabla_0$ on ${\mathcal O}_X$ such that $({\mathcal L},\nabla)$ is isomorphic to 
$({\mathcal O}_X,\nabla_0)$.  So $\nabla_0$ is a regular connection. On the other hand
there is a global form $\omega$ on $X$, such that $\nabla_0(f)=df+\omega f$. If $\nabla_0$ is regular,
one can choose the form $\omega$ as a global rational form on $\bar{X}$ in 
$H^0(\bar{X},\Omega^1_{\bar{X}}(log\,D))$. Then the lower row is exact at $Pic_{cr}(X)$.

Now, let us check the exactness at $Pic(X)$. Let $\bar{\mathcal U}=\bar{U}_i$ be an affine covering  of $\bar{X}$ as
in the proof of Lemma \ref{exact}, such that $(U_i)$ is a covering of $X$ and 
$({\mathcal L}|{U}_i,\nabla|{U}_i=\nabla_i)$
is isomorphic to $({\mathcal O}_X|{U}_i,\tilde{\nabla}_i)$, where:
$$\tilde{\nabla}_i(f)=df + \omega_i f$$
with a rational differential form $\omega_i$ defined on $\bar{U}_i$ with poles contained in $D$. On 
this covering $(U_i)$ of $X$, the invertible sheaf ${\mathcal L}$ defines the cocycle $(g_{ij})$ and 
its image in $H^1(\bar{X},\Omega^1_{\bar{X}}(log\,D))$ is the cocycle  
$\frac{d\hat{g}_{ij}}{\hat{g}_{ij}}$ defined by the rational functions on the covering $(\bar{U}_i)$
which extend $(g_{ij})$ and, again:
$$\frac{d\hat{g}_{ij}}{\hat{g}_{ij}}=\omega_j-\omega_i.$$
It remains to prove the exactness at $Pic(X)$.

If the image of the class of ${\mathcal L}$ in $H^1(\bar{X},\Omega^1_{\bar{X}}(log\,D))$
is trivial, we have:
$$\frac{d\hat{g}_{ij}}{\hat{g}_{ij}}=\omega_j-\omega_i.$$
where $\hat{g}_{ij}$ is a rational function which extends $g_{ij}$ to $\bar{X}$ and $\omega_i$
is a logarithmic differential form along $D$ on $\bar{U}_i$. The invertible sheaf
${\mathcal L}$ is endowed with a regular connection $\nabla$ locally defined on $U_i$ by:
$$\tilde{\nabla}_i(f)=df + \omega_i f.$$
This ends the proof of Lemma \ref{log}.\\

Then, we have:

\begin{lemma} $H^0(\bar{X},{^c\Omega}^1_{\bar{X}}(log\,D))=
H^0(\bar{X},\Omega^1_{\bar{X}}(log\,D))$
\end{lemma}

\noindent{\bf Proof.} We know that the spectral sequence $E_1^{pq}=
H^q(\bar{X},\Omega^p_{\bar{X}}(log\,D))\to H^{p+q}(X^{an};\mathbb{C})$ 
degenerates at $E_1$ 
(see \cite{[D2]} Corollaire 3.2.13 page 38), so the mapping:
$$H^0(\bar{X},\Omega^1_{\bar{X}}(log\,D))\stackrel{d}{\rightarrow}
H^0(\bar{X},\Omega^2_{\bar{X}}(log\,D))$$ 
is the zero map which precisely means that the forms in 
$H^0(\bar{X},\Omega^1_{\bar{X}}(log\,D))$
are closed as stated in the lemma.\\
This proves the Lemma.\\

Now let $I$ be the image of $Pic\,X$ in $H^1(\bar{X},{^c\Omega}^1_{\bar{X}}(log\,D))$. 
We shall show that the homomorphism
$I\to H^1(\bar{X},\Omega^1_{\bar{X}}(log\,D))$ is injective.\\

Consider the complex of sheaves $F^1\Omega^*_{\bar{X}}(log D)$ which is the first
term of the filtration b\^ete of $\Omega^1_{\bar{X}}(log D)$. There is a natural arrow
from the simple complex with $^c\Omega^1_{\bar{X}}(log D)$ in degree $0$ into the complex
$F^1\Omega^*_{\bar{X}}(log D)$ shifted by $-1$ which gives an exact sequence of complexes:
$$0\rightarrow {^c\Omega^1_{\bar{X}}}(log D)\rightarrow F^1\Omega^*_{\bar{X}}(log D)(-1)
\rightarrow {\mathcal C}^{\cdot}\rightarrow 0.$$
where ${\mathcal C}^{\cdot}$ is the following non-negative complex:
$$\Omega^1_{\bar{X}}(log\,D)/{^c\Omega}^1_{\bar{X}}(log\,D)\to
 \Omega^2_{\bar{X}}(log\,D)\to\ldots$$
We notice that the mapping $H^1(\bar{X},{^c\Omega}^1_{\bar{X}}(log\,D))\to
 \mathbb{H}^1(\bar{X},F^1\Omega^*_{\bar{X}}(log D)(-1))$ is injective, since the $0$-cohomology
 of ${\mathcal C}^\cdot$ is obviously $0$, so the hypercohomology
 ${\mathbb H}^0(\bar{X},{\mathcal C}^{\cdot})=0$.
 
Furthermore, by the property of the filtration b\^ete: 
$$\mathbb{H}^2(\bar{X},F^1\Omega^*_{\bar{X}}(log D))=F^1H^2(X^{an};\mathbb{C})=F^1.$$ 

We saw that $H^1(\bar{X},{^c\Omega}^1_{\bar{X}}(log\,D))$ injects in
 $\mathbb{H}^1(\bar{X},F^1\Omega^*_{\bar{X}}(log D)(-1))=
\mathbb{H}^2(\bar{X},F^1\Omega^*_{\bar{X}}(log D))$ therefore $I$ injects in
$F^1H^2(X^{an};\mathbb{C})=F^1$.  \\

Then the image of $I$ in $H^2(X^{an};\mathbb{C})$ is in $F^1$.\\

We have a commutative diagram (beware the upper line is not exact):
$$\begin{array}{ccccc}
Pic(\bar{X}) &\to &Pic(X) &\to &H^1(\bar{X},{^c\Omega}^1_{\bar{X}}(log\,D)) \\
\downarrow &&\downarrow &\swarrow&\\
H^2(\bar{X}^{an},{\mathbb C}) &\to &H^2(X^{an},{\mathbb C})&&
\end{array}$$

By e.g. \cite{[HL1]} (p. 75), the map $Pic(\bar{X}) \to Pic(X)$ is surjective, so $I$ is also
the image of $Pic(\bar{X})$ in $H^1(\bar{X},{^c\Omega}^1_{\bar{X}}(log\,D))$.
Since the cohomology of the smooth compactification $\bar{X}$ of $X$ has a pure Hodge structure, we can conclude that the image of $I$ in $H^2(X^{an};\mathbb{C})$ is contained in $W_2$. Therefore this 
image of $I$ is in the subspace $F^1\cap\overline{F}^1\cap W_2$ of 
$H^2(X^{an};\mathbb{C})$. This coincides with $Gr_F^1W_2$, which is a subset of $Gr_F^1=H^1(\bar{X},\Omega^1_{\bar{X}}(log\,D))$.\\
So we conclude that $I\to H^1(\bar{X},\Omega^1_{\bar{X}}(log\,D))$ is injective.

Now we apply Five Lemma, as in the proof of Theorem \ref{kahler}. So $Pic_{cir}X=Pic_{cr}X$.

Finally suppose that a regular connection on $\mathcal L$ is given. As we have proved there is an invertible $\mathcal{O}_X$-module ${\mathcal L}'$ which is isomorphic to $\mathcal{L}$ and which admits a regular integrable connection. The isomorphism induces a regular integrable connection on $\mathcal L$. Both connections on $\mathcal L$ differ by an element of $H^0(\bar{X},\Omega^1_{\bar{X}}(log\,D))$. By the degeneration of the Hodge spectral sequence we know that this form is closed, so both connections have the same curvature. Therefore our original regular connection must be integrable, too.

This ends the proof of Theorem \ref{int}.

\vskip.1in

\noindent {\bf Remark:} Since $Pic_{cir}\,X=Pic_{ci}\,X$ we have $Pic_{cir}\,X\simeq \mathbb{H}^1(X,\tilde{\mathcal T}^\cdot)$, where $\tilde{\mathcal T}^\cdot$ is the complex $j_*{\mathcal O}^*_X\to\Omega^1_{\overline{X}}(log\,D)\to 0\to \ldots$.

  \section{Some examples} 
 
\vskip.1in 
\noindent In the following example we only consider complex algebraic varieties.
  
 \subsection{} For the complex projective line, the invertible sheaf ${\mathcal O}(k)$ has
no connection whenever $k\neq 0$. One knows that $Pic({\mathbb P}^1)={\mathbb Z}$. We shall see that
$Pic_{ci}(X)\simeq Pic_{c}(X)=\{0\}$. 

In fact, as we have proved in the section 2, for any compact connected complex K\"ahler manifold $X$ 
(in particular any complex projective variety without singularities) we have:
$$Pic^{an}_{ci}(X)\simeq Pic^{an}_{c}(X)\simeq H^1(X,{\mathbb C}^*).$$
For the complex line ${\mathbb P}^1$ the cohomology $ H^1(X^{an},{\mathbb C}^*)=0$. By GAGA 
(see \cite{S}, \cite{M}) we have 
$Pic_{ci}(X)\simeq Pic^{an}_{ci}(X^{an})$ and $Pic_{c}(X)\simeq Pic^{an}_{c}(X^{an})$. 

\subsection{} \label{complicated} We give an example of an invertible $\mathcal{O}_X$-module which has a connection but no integrable connection.

Let $\bar{X}:=\{z_0z_1-z_2z_3=0\}\subset{\mathbb P}^3$. Notice that $\bar{X}$ 
is a complex surface isomorphic to ${\mathbb P}^1\times{\mathbb P}^1$.

Let $D:=\bar{X}\cap\{z_0+z_1+z_2-z_3=0\}$. Let $X:=\bar{X}\setminus D$.

One verifies that $D$ is a smooth hypersurface of $\bar{X}$. Using Lefschetz Theorem on hyperplane sections, one shows that $D$ is connected. In fact, $D$ is a non-singular projective plane 
curve of degree $2$. So $D\simeq {\mathbb P}^1$. Then $H^1(D^{an};{\mathbb Z})=0$.

By \cite{[HL1]} (p. 75) we have a commutative diagram whose lines are exact:
$$\begin{array}{ccccccc}
{\mathbb Z}&\to&Pic\,\bar{X}&\to&Pic\,{X}&\to&0\\
\downarrow\simeq&&\downarrow&&\downarrow&&\\
H^2(\bar{X}^{an},X^{an};{\mathbb Z})&\to&H^2(\bar{X}^{an};{\mathbb Z})&\to&Im\, \phi&\to &0
\end{array}$$
where $\phi:H^2(\bar{X}^{an};{\mathbb Z})\to H^2({X}^{an};{\mathbb Z})$.

We have (see \cite{[H]} Chap. III Exercise 12.6, p. 292):
$$Pic\,\bar{X}\simeq Pic\,{\mathbb P}^1\times Pic\,{\mathbb P}^1\simeq 
{\mathbb Z}\times{\mathbb Z}.$$
According to K\"unneth formula, we have:
$$H^2(\bar{X}^{an};{\mathbb Z})\simeq{\mathbb Z}\oplus{\mathbb Z}.$$
One verifies that the middle vertical arrow in the diagram above given by 
the first Chern class is an isomorphism: one has to compute 
$c_1(p^*_i({\mathcal O}_{{\mathbb P}^1}(n))$,
$i=1,2$, where $p_1$ and $p_2$ are the projections of $\bar{X}$ onto ${\mathbb P}^1$.

By the Five Lemma, the last vertical arrow is an isomorphism. 

Moreover the lower line of the diagram gives an exact sequence:
$${\mathbb Z}\to{\mathbb Z}\oplus {\mathbb Z}\to Im\,\phi$$
because $H^2(\bar{X}^{an},X^{an};{\mathbb Z})\simeq H_2(D^{an};{\mathbb Z})$ 
by Lefschetz duality
and:
$$H_2(D^{an};{\mathbb Z})\simeq {\mathbb Z}$$
because $D\simeq {\mathbb P}^1$.

Therefore, there is an element $c\in Im\,\phi$ which is not a torsion element.

The surjectivity of the third vertical arrow gives that there is an invertible sheaf ${\mathcal L}$ on $X$ such that $c_1({\mathcal L})=c$.

Since $X$ is affine, we have:
$$H^1(X,\Omega^1_X)=0.$$
According to Lemma \ref{special2} there is a connection on the sheaf
${\mathcal L}$. But according to Lemma
\ref{chern}, there is no integrable connection on ${\mathcal L}$.

In fact we can be more precise in this case. The lower line of the diagram of \cite{[HL1]} (p. 75) is:
$$H^2(\bar{X}^{an},X^{an};{\mathbb Z})\to H^2(\bar{X}^{an};{\mathbb Z})\to 
H^2(X^{an};{\mathbb Z})\to H^3(\bar{X}^{an},X^{an}; {\mathbb Z}).$$
By Lefschetz duality, we have:
$$H^3(\bar{X}^{an},X^{an}; {\mathbb Z})\simeq H_1(D^{an};{\mathbb Z})$$
which is $0$.
Therefore, the following sequence is exact:
$$H^2(\bar{X}^{an},X^{an};{\mathbb Z})\to H^2(\bar{X}^{an};{\mathbb Z})\to H^2(X^{an};{\mathbb Z})\to  0.$$
However $X$ is isomorphic to the set of points $(z_0:\ldots :z_3)$ of ${\mathbb P}^3$ such that: 
$$(z_0-z_1-z_2+z_3)z_1-z_2z_3=0$$ 
and $z_0\neq 0$. Therefore $X$ is isomorphic to the the variety of
points $(z_1,z_2,z_3)$ of ${\mathbb C}^3$ where $(1-z_1-z_2+z_3)z_1-z_2z_3=0$ or
isomorphic to the variety of points $(z_0,z_2,z_2,z_3)$ of ${\mathbb C}^4$ where:
$$(z_0-z_1-z_2+z_3)z_1-z_2z_3=0 \hbox{ and }z_0=1.$$
Then, one sees that $X$ is homeomorphic to the Milnor fiber of $z_0$ restricted to the hypersurface
$(z_0-z_1-z_2+z_3)z_1-z_2z_3=0$ at the point $0$. Therefore it is homemorphic to a bouquet of one sphere, i.e. one real sphere:
$$H^2(X^{an};{\mathbb Z})={\mathbb Z}.$$
Then, we have:
$$\begin{array}{ccccccc}
H^2(\bar{X}^{an},X^{an};{\mathbb Z})&\to &H^2(\bar{X}^{an};{\mathbb Z})&\to 
&H^2(X^{an};{\mathbb Z})&\to &0\\
\downarrow\simeq &&\downarrow\simeq &&\downarrow\simeq && \\
{\mathbb Z}&\to &{\mathbb Z}\oplus{\mathbb Z}&\to &{\mathbb Z}&\to 0&
\end{array}$$
 Therefore the map $H^2(\bar{X}^{an},X^{an};{\mathbb Z})\to H^2(\bar{X}^{an};{\mathbb Z})$
 is injective.

Then we obtain the $Pic(X)={\mathbb Z}$ and $Pic(\bar{X})={\mathbb Z}\oplus{\mathbb Z}$.

Since $X$ is affine, $H^1(X,\Omega^1_X)=0$. The lower line of the diagram in Theorem \ref{alg}
gives that $Pic_c(X)\to Pic(X)$ is surjective. Using Lemma \ref{chern}, the map:
$$Pic_{ci}(X)\to Pic(X)$$
is not surjective, since there are invertible sheaves which have a first Chern class which is 
not a torsion element.

By the way, we can observe that $Pic^{an} (X^{an})$ 
is isomorphic to $H^2(X^{an};{\mathbb Z})$ because
$X^{an}$ is a Stein space and $H^2(X^{an},{\mathcal O}_X)=0$. So:
$$Pic^{an} (X^{an}) \simeq{\mathbb Z} .$$

\subsection{} Notice that it is easier to find that there are connections which 
are not integrable or regular. One may consider
$X={\mathbb C}^2$. In this case both $Pic^{an}(X^{an})$ and $Pic(X)$
are trivial.

A connection on ${\mathcal O}_X$ (resp.  ${\mathcal O}_{X^{an}}$) 
is given by a global algebraic (resp. analytic) differential form
$\omega$:
$$\nabla (f)=df + \omega f.$$

If one considers $\omega=dz_1$, the corresponding connection is integrable but not regular.

If $\omega=z_1dz_2$, the corresponding connection is not integrable because the form is not closed.

We can compute $Pic_{c}(X)$ and $Pic_{ci}(X)$ by using the diagram of
Theorem \ref{alg}. Then:
$$Pic_{c}(X)\simeq H^0(X,\Omega^1_X)$$
because for $X={\mathbb C}^2$, the map 
$H^0(X,{\mathbb C}^*_{X})\to H^0(X,{\mathcal O}^*_X)$ is
an isomorphism.

Similarly, we have:
$$Pic_{ci}(X)\simeq H^0(X, {^c\Omega}^1_X).$$

In the analytic case, we know that:
$$Pic^{an}_{ci}(X^{an})\simeq H^1(X^{an},{\mathbb C}^*),$$
so it is trivial.

For $Pic^{an}_c(X^{an})$ the exact sequence of \ref{CI} gives that $Pic^{an}_c(X^{an})$ is
isomorphic to $H^0(X^{an},d\Omega^1_{X^{an}})$. The elements of 
$Pic^{an}_c(X^{an})$ are given by their curvature.

\subsection{} Consider the algebraic variety $X={\mathbb C}^*\times {\mathbb C}^*$.\\

Notice that for this variety $Pic(X)=0$, because $X={\mathbb C}^2\setminus Z$ where
$Z$ is the closed algebraic subspace given by the union of the lines ${\mathbb C}\times\{0\}$
and $\{0\}\times {\mathbb C}$, then by the Proposition 6.5 in Chapter II of \cite{[H]} p. 133, we have a surjection:
$$Pic({\mathbb C}^2)\rightarrow Pic(X).$$
Then, $Pic(X)=0$.

On the other hand $Pic^{an}(X^{an})\simeq H^2(X^{an};{\mathbb Z})={\mathbb Z}$ because 
$X^{an}$ is a Stein space and using the exact exponential sequence.

Therefore, there are invertible $\mathcal{O}_{X^{an}}$-modules for which the complex first Chern class
is $\neq 0$. According to Lemma \ref{chern} these sheaves do not have integrable connections.
However, by Theorem \ref{SE} they have a connection because $H^1(X^{an},\Omega^1_{X^{an}})=0$.
But these do not come from an algebraic invertible sheaf, because the latter ones are trivial.

\subsection{} Put $X:=\mathbb{C}^2\setminus\{0\}$. Note that $X^{an}$ is simply connected. \\
On the other hand, $H^1(X^{an},{\mathcal O}_{X^{an}})\neq 0$: Let $\mathcal U$ be the open Stein covering by $U_1=\mathbb{C}\times {\mathbb C}^*$, $U_2={\mathbb C}^*\times \mathbb{C}$. Then $H^1(X^{an},{\mathcal O}_{X^{an}})$ is the cokernel of :
$$H^0(U_1,{\mathcal O}_{U_1})\oplus H^0(U_2,{\mathcal O}_{U_1})\to H^0(U_1\cap U_2,{\mathcal O}_{U_1\cap U_2})$$ 
$$(a,b)\mapsto r_1(a)-r_2(b)$$ 
where $r_1,r_2$ are restrictions, so it corresponds to all globally convergent Laurent series in two variables with negative exponents.

The exact sequence 
$$0=H^1(X;\mathbb{Z})\to H^1(X^{an},{\mathcal O}_{X^{an}})\to Pic_0(X^{an})\to 0$$ 
shows that $Pic_0(X^{an})\neq 0$. On the other hand, $Pic(X)=Pic(\mathbb{C}^2)=0$. So there are invertible $\mathcal{O}_{X^{an}}$-modules with Chern class $0$ which are not algebraizable. These cannot admit a connection: The composition $H^1(X^{an},\mathcal{O}_{X^{an}})\to Pic(X^{an})\to H^1(X^{an},\Omega^1_{X^{an}})$ is given by $(f_{ij})\mapsto (2\pi idf_{ij})$, so $b(\mathcal{L})\neq 0$ if $(f_{ij})$ does not represent the trivial element. 

In particular, we cannot improve Lemma \ref{chern} in general. On the other hand, cf. Lemma \ref{integrable} and Corollary \ref{tors1}.

\subsection{}
 Let $X$ be a complex algebraic variety, 
$\mathcal L$ an invertible $\mathcal{O}_X$-module, $\nabla$ a connection on $\mathcal L$.\\
Then we have:
$$\nabla\hbox{ regular integrable }\Rightarrow \nabla\hbox{ integrable}$$
This implication is not invertible, as shown by the example 
$X=\mathbb{C}^2, {\mathcal L}={\mathcal O}_X$ (see above 4.3). \\
Note that $\nabla$ is integrable if and only if $\nabla^{an}$ is integrable.\\
In fact, we can consider the existence of connections on $\mathcal L$ (resp. ${\mathcal L}^{an}$):
$$\begin{array}{ccccc}
\exists\,\nabla \hbox{ regular integrable}&\Leftrightarrow &\exists\, \nabla \hbox{ integrable}&\Rightarrow &\exists\, \nabla\\
&&\Updownarrow&&\Downarrow\\
&&\exists\,\nabla \hbox{ analytic integrable}&\Rightarrow &\exists\,\nabla\hbox{ analytic}
\end{array}$$
For the left upper and the middle vertical equivalence see Corollary \ref{tors1}.\\
Note that there may be no connection at all on $\mathcal{L}$ or $\mathcal{L}^{an}$, as shown by the example $X=\mathbb{P}_1, {\mathcal L}={\mathcal O}(k), k\neq 0$.\\
The right horizontal arrows are not invertible, as shown by the complicated example \ref{complicated}.\\
The right vertical arrow is not invertible if the answer to the following question is positive:\\
Let $X$ be the Serre example of a non-singular algebraic surface which is not affine but the corresponding complex analytic manifold is Stein (see \cite{[Ha]} p. 232 Example 3.2). Is there an invertible $\mathcal{O}_X$-module $\mathcal L$ on $X$ such that its image in $H^1(X,\Omega^1_X)$ does not vanish?
(Note that $X$ is not affine, so it is possible that $H^1(X,\Omega^1_X)\neq 0$). Then, $\mathcal L$ does not admit a connection. \\
On the other hand, $X^{an}$ is Stein, so $H^1(X^{an},\Omega^1_{X^ {an}})=0$, which implies that there is a connection on ${\mathcal L}^{an}$.

\end{document}